\documentclass[a4paper,11pt]{article}

 \usepackage{a4wide}
 
\usepackage{amsmath}
\usepackage{amsfonts}
\usepackage{amsthm}

\usepackage{graphicx}
\usepackage{color}

\def\R{\mathbb{R}}
\def\C{\mathbb{C}}
\def\Z{\mathbb{Z}}
\def\N{\mathbb{N}}
\def\pa{\partial}
\def\CO{{\mathcal O}}
\def\CP{{\mathcal P}}
\def\la{\langle}
\def\ra{\rangle}

\def\be{\begin{equation}}
\def\ee{\end{equation}}

\def\Op{\mathop{{\rm Op}^W_{h}}\nolimits}

\renewcommand{\Re}{\text{{\rm Re}\;}}
\renewcommand{\Im}{\text{{\rm Im}\;}}

\newcommand{\supp}{\text{{\rm Supp}}\hskip 1pt}

\newcommand{\dist}{\text{\rm dist}\hskip 1pt}

\newtheorem{thm}{Theorem}[section]
\newtheorem{lem}[thm]{Lemma}
\newtheorem{prop}[thm]{Proposition}

\newtheorem{rem}[thm]{Remark}

\newtheorem{definition}[thm]{Definition}

\numberwithin{equation}{section}

\title{Resonances for non-analytic potentials}
\author{Andr\'e Martinez${}^1$, Thierry Ramond${}^2$, and  Johannes Sj\"ostrand${}^3$}

\begin{document}

\maketitle 
\addtocounter{footnote}{1}
\footnotetext{Universit\`a di Bologna, Dipartimento di
Matematica, Piazza di Porta San Donato 5, 40127 Bologna,
Italy. Partly supported by Universit\`a di Bologna, Funds
for Selected Research Topics and Founds for Agreements with
Foreign Universities}
\addtocounter{footnote}{1}
\footnotetext{D\'epartement de Math\'ematiques, Universit\'e Paris-Sud 11, UMR CNRS 8628, FR-91405 Orsay, France.}
\addtocounter{footnote}{1}
\footnotetext{CMLS, Ecole Polytechnique, UMR CNRS  7640, FR-91128 Palaiseau Cedex, France.}

%%%%%%%%%%%%%%%% ABSTRACT %%%%%%%%%%%%%%%%%%%%%%%%%%%%%
\begin{abstract} 
We consider semiclassical Schr\"odinger operators on $\R^n$, with $C^\infty$ potentials decaying polynomially at infinity. The usual theories of resonances do not apply in such a non-analytic framework. Here, under some additional conditions, we show that resonances are invariantly defined up to any power of their imaginary part. The theory is based on resolvent estimates for  families of approximating distorted operators with potentials that are holomorphic in narrow complex sectors around $\R^n$.
\end{abstract}

%%%%%%%%%%%%%%%%%  SECTION 1  %%%%%%%%%%%%%%%%%%%%%%%%%

 \section{Introduction} \label{sec-intro}

In physics, the notion of quantum resonance has appeared at  the begining of quantum mechanics. Its introduction was motivated by the behavior of various quantities related to scattering experiments, such as the scattering cross-section. At certain energies, these quantities present peaks (nowaday called Breit-Wigner peaks),  which were modelized  by  a Lorentzian shaped function 
$$
w_{a,b} : \lambda \mapsto ((\lambda-a)^2 +b^2)^{-1}.
$$
The real numbers $a$ and $b$ stand for the location of the maximum of the peak and its height. Of course for $\rho = a-ib\in \C$, one has 
$$
 w_{a,b} (\lambda )= \frac{1}{\vert \lambda -\rho\vert^2},
$$
and the complex number $\rho$ was called  a resonance. Such complex values for energies had also appeared for example in the work  \cite{ga} by Gamow, to explain $\alpha$-radioactivity, and were associated to the existence of some decaying state of energy $a=\Re \rho$ and lifetime $1/b=1/|\Im \rho|$.

However,  these complex numbers are not defined in a completely exact way, in the sense that the peaks do not perceivably change if these numbers are modified by a quantity much smaller than their imaginary part. Indeed, a  straightforward computation shows that the relative difference between such two peaks $w_{a,b}$ and $w_{a',b'}$  verifies,
$$
\sup_{\lambda\in \R}\left\vert\frac{ w_{a,b} (\lambda )- w_{a',b'} (\lambda )}{w_{a',b'}(\lambda)}\right\vert \leq 2\left\vert\frac{\rho -\rho'}{\Im\rho}\right\vert+ \left\vert\frac{\rho -\rho'}{\Im\rho}\right\vert^2
$$
where we have also set $\rho'=a'-ib'$. 
As a consequence, the two peaks become undistinguishable if $\vert \rho -\rho'\vert << \vert\Im\rho\vert$, that is, there is no physical relevance to associate the resonance $\rho =a-ib$ to $w_{a,b}$ rather than any other $\rho'$ verifying $
\vert \rho -\rho'\vert << \vert\Im\rho\vert$. Notice also that  the more the resonance is far from the real line, the more  irrelevant this precision becomes.

On the mathematical side, 
the more recent  theory of  resonances  for Schr\"odinger operators has permitted to give a rigorous  framework and to obtain very precise results, in particular on the location of resonances in relation with the geometry of the underlying classical flow. However, it is based on the notion of complex scaling, in more and more sophisticated versions (see, e.g., 
 \cite{AgCo, BaCo, Sim, Hu, Sig, Cy, na1, na2, HeSj})
 that all require analyticity assumptions on the potential (or its Fourier transform). 

There is a small number of  works about the definition of resonances for non-analytic potentials, as  e.g. \cite{Or,GeSi, SoWe, JeNe, CaMaRa}.  In \cite{Or,GeSi, SoWe, JeNe}, the point of view is quite different from ours, while in
\cite{CaMaRa}, the definition is based on the use of an almost-analytic
extension of the potential and seems to strongly depend both on the
choice of this extension and on the complex distortion.

 Here our purpose is to  give a definition that  fulfills both the mathematical requirement of being invariant with respect to the choices one has to make, and the physical requirement of being more accurate as the resonance become closer to the real (or, equivalently, as the Breit-Wigner peak becomes narrower). Dropping the physically irrelevant precision for the definition of resonances, we can also drop the spurious assumption on the analyticity of the potential.
 
 More precisely, we associate to a Schr\"odinger operator $P$  a discrete set $\Lambda\subset \C$ with certain properties, such that, for any other set $\Lambda'$ with the same properties, there exists a bijection $B : \Lambda' \rightarrow \Lambda$ with $B(\rho) -\rho =\CO (\vert \Im\rho\vert^\infty)$ uniformly. The set of resonances of $P$ is the corresponding equivalence class of $\Lambda$. Of course, when the potential is dilation analytic at infinity, we recover the usual set of resonances up to the same error $\CO (\vert \Im\rho\vert^\infty)$.

 The properties characterizing $\Lambda$ basically involve the resonances of a (essentially arbitrary) family of dilation-analytic operators $(P^\mu)_{0<\mu\leq\mu_0)}$, such that,
 \begin{eqnarray*}
 && P^\mu \mbox{ is dilation-analytic in a complex sector of angle } \mu \mbox{ around } \R^n;\\
 && \Vert P^\mu -P\Vert =\CO(\mu^\infty ) \mbox{ uniformly as } \mu\rightarrow 0_+,
 \end{eqnarray*}
and the constructive proof of the existence of the set $\Lambda$ mainly consists in studying  such a family and, in particular, in obtaining resolvent estimates uniform in $\mu$.

In this paper, we address the case of  an isolated cluster  of resonances with a bounded (with respect to $h$)  cardinality.
We hope to treat  the general case elsewhere, as well as  to give a detailed description of the quantum evolution $e^{itP/h}=e^{itP^\mu/h}+\CO(|t|h^{-1}\mu^\infty)$ in terms of the resonances in $\Lambda$.

The paper is organized as follows. We give our assumptions and state our main results  in Section \ref{sec-result}. Then, in Section 3, we give two paradigmatic situations where our constructions apply: the non-traping case and the shape resonances case. In section 4 we present a suitable notion of  analytic approximation of a $C^\infty$ function through which we define the operator $P^\mu$. In Section \ref{sec-dist} we show that  a properly defined analytic distorted operator $P^\mu_{\theta}$ of the latter verifies a nice resolvent estimate in the upper half complex plane even very near to the real axis.
The sections \ref{sec-th21},\ref{sec-red} and \ref{proofTh3} are devoted to the proof of Theorem \ref{Th1}, Theorem \ref{Th2} and Theorem \ref{Th3} respectively. We construct the set of resonances $\Lambda$, and prove Theorem \ref{Th4} in Section \ref{cons}. In the last Section \ref{sec-exemp}, we prove our statements concerning the shape resonances. Eventually, we have placed in Appendix \ref{sec-app} the proofs of two technical lemmas.

%%%%%%%%%%%%%%%%%%%%%%%  SECTION 2  %%%%%%%%%%%%%%%%%%%%%%%%%%%%%%%
\section{Notations and Main Results}\label{sec-result}

We consider the semiclassical
 Schr\"odinger operator,
%\begin{equation}\label{eq-result-1}

$$
P= -h^2\Delta + V,
$$
%\end{equation}
where $V=V(x)$ is a real smooth function of $x\in\R^n$, such that,
\be
\label{decV}
\pa^\alpha V(x) =\CO (\la x\ra^{-\nu-|\alpha|}),
\ee
for some $\nu >0$ and for all $\alpha\in\Z_+^n$.
 We also  fix $\widetilde\nu\in (0,\nu)$ once for all, and, for any $\mu >0$ small enough, we denote by $V^{\mu}$ a $|x|$-analytic $(\mu, \widetilde\nu)$-approximation of $V$ in the sense of Section \ref{sec-prelim}. In particular,  $V^\mu$ is analytic with respect to $r=|x|$ in $\{r\geq 1\}$, it can  be extended into a holomorphic function of $r$ in the sector $\Sigma:=\{ \Re r\geq 1\, ,\, |\Im r|\leq 2\mu\Re r\}$, and it verifies,
\be
\label{compVVhol1}
 V^\mu(x) -  V(x) = \CO ( \mu^\infty\la x\ra^{-\widetilde\nu }),
\ee
uniformly on $\R^n$. (See Section \ref{sec-prelim} for more properties of $V^\mu$.)

Then, for any $\theta \in (0,\mu]$, the operator,
\be
P^\mu:= -h^2\Delta + V^\mu,
\ee
can be distorded analytically into,
\be
\label{distP}
P^\mu_{\theta}:= U_{\theta}P^\mu U_{\theta}^{-1},
\ee
where $U_\theta$ is any transformation of the type,
\be
\label{Utheta}
U_\theta \varphi (x) := \varphi (x+i\theta A(x)),
\ee
with $A(x):=a(|x|)x$, $a\in C^\infty (\R_+)$, $a=0$ near 0, $0\leq a\leq 1$ everywhere, $a(|x|) =1$ for $|x|$ large enough.
The essential spectrum of $P^\mu_{\theta}$ is $e^{-2i\theta}\R$, and
its discrete spectrum $\sigma_{disc}(P^\mu_{\theta})$ is included in the lower half-plane and does not depend on the choice of the function $a$. Moreover, it does not depend on $\theta$,  in the sense that for any $\theta_0\in (0,\mu]$, and  any $\theta\in [\theta_0,\mu]$ , one has,
\begin{eqnarray*}
&& \sigma_{disc}(P^\mu_{\theta})\cap \Sigma_{\theta_0}=\sigma_{disc}(P^\mu_{\theta_0})\cap \Sigma_{\theta_0},
\end{eqnarray*}
where we have set $\Sigma_{\theta_0}:=\{z\in \C\; ; \; -2\theta_0< \arg z\leq 0\}$ (observe that one also has $\sigma_{disc}(P^\mu_{\theta}) =\sigma_{disc}(\widetilde U_{\theta}P^\mu\widetilde U_{\theta}^{-1})$, where $\widetilde U_\theta \varphi (x) := \sqrt{\det (Id + i\theta\hskip 1pt {}^t dA(x))} \varphi (x+i\theta A(x))$ is an analytic distorsion more widely used in the literature).

We denote by,
$$
\Gamma (P^\mu) := \sigma_{disc}(P^\mu_{\mu})\cap \Sigma_{\mu},
$$
the set of resonances of $P^{\mu}$ counted with their multiplicity. In what follows, we also  use the following notation: If $E$ and $E'$ are two $h$-dependent subsets of $\C$, and $\alpha =\alpha (h)$ is a $h$-dependent positive quantity that tends to 0 as $h$ tends to $0_+$, we write,
$$
E' =E +{\mathcal O}(\alpha),
$$
when there exists a constant $C>0$ (uniform with respect to all other parameters) and a bijection
$$
b : E' \to E,
$$
such that,
$$
\vert b(\lambda) -\lambda\vert \leq C\alpha
$$
for all $h>0$ small enough.

Now, we fix some energy level $\lambda_0>0$, and a constant $\delta >0$. For any $h$-dependent numbers $\widetilde\mu(h), \mu(h)$, and any  $h$-dependent bounded intervals $I(h), J(h)$, verifying,
\begin{eqnarray}
\label{condmu}
&& 0 <\widetilde\mu(h) \leq \mu(h)\leq h^{\delta};\\
\label{condint}
&& I(h)\subset J(h)\; ; \; \mbox{ diam}(J\cup \{\lambda_{0}\})\leq h^{\delta},
\end{eqnarray} 
we consider the following property:

\bigskip
${\CP}(\widetilde\mu,\mu ; I,J):$
\qquad $\left\{
\begin{minipage}{12cm}
$\begin{array}{l}
\Re(\Gamma(P^\mu)\cap (J-i[0,\lambda_{0}\widetilde{\mu}])) \subset I;\\[6pt]
\# (\Gamma(P^\mu)\cap (J-i[0,\lambda_{0}\widetilde{\mu}]))\leq \delta^{-1};\\[6pt]
\dist( I,\R\setminus J)\geq h^{-\delta} \omega_{h}(\widetilde{\mu}),
\end{array}
$
\end{minipage}
\right.
$

\begin{figure}[h] %
\begin{center}
\begin{picture}(0,0)%
\includegraphics{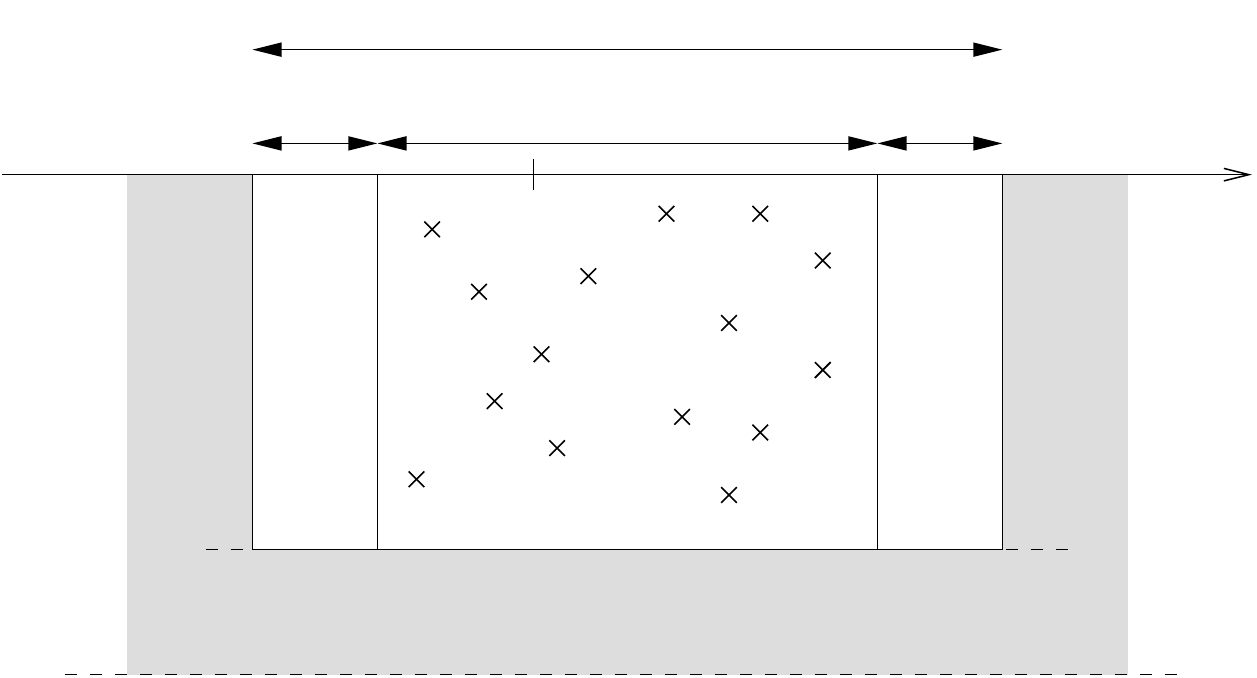}%
\end{picture}%
\setlength{\unitlength}{3947sp}%
\begingroup\makeatletter\ifx\SetFigFont\undefined%
\gdef\SetFigFont#1#2#3#4#5{%
  \reset@font\fontsize{#1}{#2pt}%
  \fontfamily{#3}\fontseries{#4}\fontshape{#5}%
  \selectfont}%
\fi\endgroup%
\begin{picture}(6024,3246)(1789,-2773)
\put(2926,-61){\makebox(0,0)[lb]{\smash{{\SetFigFont{12}{14.4}{\familydefault}{\mddefault}{\updefault}{\color[rgb]{0,0,0}$h^{-\delta}\omega_h(\tilde\mu)$}%
}}}}
\put(4276,-661){\makebox(0,0)[lb]{\smash{{\SetFigFont{12}{14.4}{\familydefault}{\mddefault}{\updefault}{\color[rgb]{0,0,0}$\lambda_0$}%
}}}}
\put(5026,-136){\makebox(0,0)[lb]{\smash{{\SetFigFont{12}{14.4}{\familydefault}{\mddefault}{\updefault}{\color[rgb]{0,0,0}$I$}%
}}}}
\put(4351,314){\makebox(0,0)[lb]{\smash{{\SetFigFont{12}{14.4}{\familydefault}{\mddefault}{\updefault}{\color[rgb]{0,0,0}$J$}%
}}}}
\put(7351,-211){\makebox(0,0)[lb]{\smash{{\SetFigFont{12}{14.4}{\familydefault}{\mddefault}{\updefault}{\color[rgb]{0,0,0}$\Re z$}%
}}}}
\put(5926,-61){\makebox(0,0)[lb]{\smash{{\SetFigFont{12}{14.4}{\familydefault}{\mddefault}{\updefault}{\color[rgb]{0,0,0}$h^{-\delta}\omega_h(\tilde\mu)$}%
}}}}
\put(6676,-2086){\makebox(0,0)[lb]{\smash{{\SetFigFont{12}{14.4}{\familydefault}{\mddefault}{\updefault}{\color[rgb]{0,0,0}$\Im z=-\lambda_0\tilde\mu$}%
}}}}
\put(6676,-2686){\makebox(0,0)[lb]{\smash{{\SetFigFont{12}{14.4}{\familydefault}{\mddefault}{\updefault}{\color[rgb]{0,0,0}$\Im z=-2\lambda_0\mu$}%
}}}}
\end{picture}%
\end{center}
\caption{ The property ${\CP}(\widetilde\mu,\mu ; I,J)$.}
 \label{fig:property}
\end{figure}

where, 
for $\theta>0$, we have set, 
$$
\omega_{h}(\theta):= \theta 
  \left (
  \ln \frac{1}{\theta} +h^{-n}(\ln \frac{1}{h})^{n+1}
   \right)^{1/2}.
$$
Notice that by (\ref{condint}), the property ${\CP}(\widetilde\mu,\mu ; I,J)$ implies  $\omega_{h}(\widetilde\mu)\leq h^{2\delta}$. 

%Then, we have,

\begin{thm}\sl \label{Th1} Suppose $\CP(\widetilde{\mu},\mu;I,J)$ holds for some $\widetilde{\mu},\mu, I$ and $J$ verifying (\ref{condmu}) -- (\ref{condint}). Then for all $\theta\in ]0,\widetilde{\mu}]$, there exists an interval
$$
J' =J+\CO(\omega_{h}(\theta)),
$$
such that,
$$
\Vert {(P^{\mu}_{\theta}-z)^{-1}}\Vert\leq C \theta^{-C}\prod_{\rho\in \Gamma(\widetilde{\mu},\mu,J)}\vert z-\rho\vert^{-1},
$$
for all $z\in J'+i[-C\theta h^{n_{1}},C\theta h^{n_{1}}]$. Here we have set $n_{1}:= n+\delta$, 
$$
\Gamma(\widetilde{\mu},\mu,J):=\Gamma(P^\mu)\cap (J-i[0,\lambda_{0}\widetilde{\mu}]),
$$
and $C>0$ is  a constant independent of $\widetilde{\mu}$, $\mu$, $\theta$, $I$ and $J$.
\end{thm}

Thanks to this result, one can compare the resonances of the operators $P^\mu$ for different values of $\mu$, as follows:

\begin{thm}\sl  \label{Th2} Let $N_{0}\geq 1$ be a constant. Suppose $\CP(\widetilde{\mu},\mu;I,J)$ holds for some $\widetilde{\mu},\mu, I$ and $J$ verifying (\ref{condmu}) -- (\ref{condint}), and that   $\widetilde{\mu}>\mu^{N_{0}}$. Then, for  any $\theta\in [\mu^{N_{0}},\widetilde{\mu}]$, there exist an interval,
$$
J' = J+\CO(\omega_{h}(\theta))
$$
and 
$\tau\in [h^{n_{1}}\theta,2h^{n_{1}}\theta]$,
such that, for any constant $N_{1}\geq 1$ and any $\mu'\in [\mu^{N_{1}}, \mu^{1/N_{1}}]$ with $\theta\leq \mu'$, one has, 
$$
\Gamma(P^{\mu'})\cap (J'-i[0,\tau]) = \Gamma(P^{\mu})\cap (J'-i[0,\tau])+\CO(\mu^{\infty}).
$$
\end{thm}

\begin{rem}\sl The only properties of $V^{\mu}$ used in the proof of this result are that $V^\mu$ is a holomorphic function of $r$ in the sector $\Sigma:=\{ \Re r\geq 1\, ,\, |\Im r|\leq 2\mu\Re r\}$,  and it verifies (\ref{compVVhol1}) and (\ref{compVVho2})  for some $\tilde\nu >0$. In particular, the proof also shows that, up to $\CO(\mu^{\infty})$, the set $\Gamma (P^\mu)$ does not depend on any particular choice of $V^\mu$.
\end{rem}

\begin{rem}\sl As we will see in the proof, the condition $\tau\in [h^{n_{1}}\theta,2h^{n_{1}}\theta]$ can actually be replaced by $\tau\in [h^{n_{1}}\theta,h^{n_{1}}\theta + (h^{n_{1}}\theta)^M]$, for any fixed $M\geq 1$.
\end{rem}

We also show that the validity of  $\CP (\widetilde{\mu},\mu; I,J)$ persists when decreasing $\widetilde{\mu}$ and $\mu$ suitably, up to a small change of $I$ and $J$.

\begin{thm}\sl  \label{Th3} Suppose $\CP(\widetilde{\mu},\mu;I,J)$ holds for some $\widetilde{\mu},\mu, I$ and $J$ verifying (\ref{condmu}) -- (\ref{condint}). Assume furthermore that there is a constant $N_{0}\geq 1$ with $\widetilde{\mu}\geq \mu^{N_{0}}$.
Then, there exist two intervals,
\begin{align*}
& I' =I+\CO(\mu^\infty);\\
& J' =J+\CO(\omega_{h}(\widetilde{\mu})),
\end{align*}
such that $\CP (h^{n_{1}}\mu',\mu';I',J')$ holds,  for any $\mu'\in (0,\widetilde{\mu}]$.
\end{thm}

Finally, the following result gives a definition of resonances for $P$, up to any power of their imaginary part.

\begin{thm}\sl  \label{Th4} Suppose $\CP(\widetilde{\mu},\mu;I,J)$ holds for some $\widetilde{\mu},\mu, I$ and $J$ verifying (\ref{condmu}) -- (\ref{condint}). Assume furthermore that there is a constant $N_{0}\geq 1$ with $\widetilde{\mu}\geq \mu^{N_{0}}$.
Then, there exist,
\begin{align*}
\mbox{an interval }\quad\, & I' = I+\CO(\mu^\infty);\\
\mbox{an interval }\quad\, & J' = J+\CO(\omega_{h}(\widetilde{\mu}));\\
\mbox{a discrete set }\, & \Lambda \subset I'-i[0,2h^{2n_{1}}\widetilde{\mu}],
\end{align*}
 such that,
 \vskip 0.5cm
\qquad($\star$)
$\left\vert
\;\begin{minipage}{13cm}
for any $\mu'\in (0,\widetilde{\mu}]$, there exist $\tau\in [h^{2n_{1}}\mu', 2h^{2n_{1}}\mu']$ with, 

$$
\Gamma(P^{\mu'})\cap (J'-i[0,\tau]) = \Lambda\cap (J'-i[0,\tau]) +\CO((\mu')^{\infty}).
$$
\end{minipage}
\right.$

\bigskip
Moreover, any other set $\widetilde{\Lambda}\subset I'-i[0,2h^{2n_{1}}\widetilde{\mu}]$ verifying ($\star$), possibly with some other choice of $V^{\mu}$, is such that there exist $\tau'\in [\frac12h^{2n_{1}}\widetilde{\mu}, h^{2n_{1}}\widetilde{\mu}]$ and a bijection,
$$
B : \Lambda\cap (J'-i[0,\tau']) \to \widetilde{\Lambda}\cap (J'-i[0,\tau']),
$$
with,
$$
B(\lambda)-\lambda=\CO(\vert \Im\lambda\vert^{\infty}).
$$
The set $\Lambda$ will be called the set of resonances of $P$ in $J'-i[0,\frac12h^{2n_{1}}\widetilde{\mu}]$. Here we adopt the convention that  real elements of $\Lambda$ are counted with a positive integer multiplicity in the natural  way (see Section \ref{cons}).
\end{thm}

\section{Two examples}\label{exemp}

Here, we describe two explicit situations where the previous results apply. 

\subsection{The non-trapping case} 

We suppose first that the energy $\lambda_0$ is non-trapping, i.e. for any $(x,\xi)\in p^{-1}(\lambda_{0})$ we have
$$
\vert \exp tH_{p}(x,\xi )\vert \to \infty \mbox{ as } \vert t\vert\to\infty,
$$
where $p(x,\xi ):= \xi^2 +V(x)$ is the principal symbol of $P$, and $H_{p}=\partial_{\xi}p\partial_{x}-\partial_{x}p\partial_{\xi}$ is the Hamilton field of $p$.

Then the result of \cite{Ma2} can be applied to $P^\mu$ with $\mu =Ch\ln (h^{-1})$  for any arbitrary constant $C>0$, and tells us that $P^\mu$ has no resonances in $[\lambda_0-2\varepsilon, \lambda_0 +2\varepsilon]-i[0,\lambda_0\mu]$ with some $\varepsilon >0$ constant. In that case, for any $\delta >0$, $\CP(h^{n_1}{\mu},\mu;I,J)$ is verified with $I=[\lambda_0-h^\delta, \lambda_0 +h^\delta]$ and $J=[\lambda_0-2h^\delta, \lambda_0 +2h^\delta]$, and the previous results tell us that $P$ has no resonances in $I-i[0, \frac12h^{3n_1}\mu]$ in the sense of Theorem \ref{Th4}.

\subsection{The shape resonances}

Now we assume instead that, in addition to (\ref{decV}), the potential $V$ presents the geometric configuration of the so-called ``point-well in an island'', as described in \cite{HeSj}. More precisely, we suppose

\bigskip

(H)\quad$\left\{\begin{minipage}{13 cm}
There exist   a connected bounded open set $\mathrm{\ddot O} \subset \R^n$, and  $x_{0}\in\mathrm{\ddot O}$, such that,
\begin{itemize}
\item
$\lambda_{0}:=V(x_{0})>0\; ;\; V> \lambda_{0} \mbox{ on } \mathrm{\ddot O}\backslash \{x_{0}\}\; ;\; \nabla^2V(x_{0}) >0\; ;\; V=\lambda_{0} \mbox{ on } \partial\mathrm{\ddot O}$;
\item Any point of $\{ (x,\xi )\in\R^{2n}\; ;\; x\in \R^n\backslash \mathrm{\ddot O}\; ,\; \xi^2 +V(x) =\lambda_{0}\}$ is non-trapping.
\end{itemize}
\end{minipage}
\right.$

\bigskip

We denote by $(e_{k})_{k\geq 1}$ the increasing sequence of (possibly multiple) eigenvalues of the harmonic oscillator $H_{0}=-\Delta +\frac12 \la V''(x_{0})x, x\ra$.
We have

\begin{thm} \sl 
\label{Th5}Assume (\ref{decV}) and (H). Then, for any $k_{0}\geq 1$ and any $\delta >0$, $\CP(\widetilde{\mu},\mu;I,J)$ holds with
$$
\mu = h^\delta\quad ; \quad \widetilde{\mu} = h^{\max(\frac{n}2,1)+1+\delta},
$$
and
$$
I=[\lambda_{0}+(e_{1}-\varepsilon )h, \lambda_{0}+(e_{k_{0}}+\varepsilon)h]\quad ; \quad J= [\lambda_{0}, \lambda_{0}+(e_{k_{0}+1}-\varepsilon)h],
$$
where $\varepsilon >0$ is any fixed number in $(0, \min (\frac{e_{1}}2, \frac{e_{k_{0}+1}-e_{k_{0}}}3)]$.
\end{thm}

Actually, we prove in Section \ref{sec-exemp} that any resonance $\rho$ of $P^{\mu}$ in $J-i[0,\widetilde{\mu}]$ is such that there exists $k\leq k_{0}$ with
$$
\Re \rho-( \lambda_{0}+e_{k}h) =\CO (h^{3/2}),
$$
and
$$
\Im\rho =\CO(e^{-2S_{1}/h}),
$$
where $S_{1}>0$ is any number less than the Agmon distance between $x_{0}$ and $\partial\mathrm{\ddot O}$. Recall that the Agmon distance is the pseudo-distance associated to the degenerate metric $(V(x)-\lambda_{0})_{+}dx^2$.

\bigskip
More generally, if $\mu'\in [e^{-\eta /h}, \mu]$ with $\eta >0$ small enough, we prove that any resonance $\rho$ of $P^{\mu'}$ in $J-i[0,\lambda_0\min(\mu', h^{2+\delta})]$, verifies
$$
\Re \rho-( \lambda_{0}+e_{k}h) =\CO (h^{3/2}),
$$
for some $k\leq k_{0}$, and
$$
\Im\rho =\CO(e^{-2(S_{0}-\eta)/h}).
$$

Applying Theorem \ref{Th4} with $\mu'=e^{-\eta /h}$ ($0<\eta<S_{0}$), we deduce that the resonances of $P$ in $[\lambda_0, \lambda_0+ Ch]-i[0, \frac12 h^{2n+\max(\frac{n}2,1) +1+3\delta}]$ satisfy the same estimates.

\section{Preliminaries}\label{sec-prelim}
 In this section, following an idea of  \cite{FuLaMa}, we define and study the notion of analytic $(\mu, \widetilde\nu)$-approximations. 
 
 \begin{definition}\sl 
 For any $\mu>0$ and $\widetilde\nu\in (0,\nu)$, we say that a real smooth function $V^\mu$ on $\R^n$ is a $|x|$-analytic $(\mu, \widetilde\nu)$-approximation of $V$, if  $V^\mu$ is analytic with respect to $r=|x|$ in $\{r\geq 1\}$, $V^\mu$ can  be extended into a holomorphic function of $r$ in the sector $\Sigma (2\mu):=\{ \Re r\geq 1\, ,\, |\Im r| < 2\mu\Re r\}$, and, for any multi-index $\alpha$, it verifies,
\be
\label{compVVho1}
\partial^\alpha ( V^\mu (x) -  V(x)) = \CO ( \mu^\infty\la x\ra^{-\widetilde\nu -|\alpha|}),
\ee
uniformly with respect to  $x\in\R^n$ and $\mu>0$ small enough, and,
\be
\label{compVVho2}
\partial^\alpha  V^\mu (x)= \CO ( \la\Re x\ra^{-\widetilde\nu -|\alpha|}),
\ee
uniformly with respect to  $x\in\Sigma (2\mu)$ and $\mu>0$ small enough.
 \end{definition}
 Then, we have,
 \begin{prop}\sl  Let $V=V(x)$ be a real smooth function of $x\in\R^n$ verifying (\ref{decV}). Then, one has,
 \begin{itemize}
 \item[(i)] For any $\mu>0$ and $\widetilde\nu\in (0,\nu)$, there exists a $|x|$-analytic $(\mu, \widetilde\nu)$-approximation of $V$;
 \item[(ii)]
 If $V^\mu$ and $W^\mu$ are two $|x|$-analytic $(\mu, \widetilde\nu)$-approximations of $V$, then, for all  $\alpha\in\N^n$, one has,
 $$
 \partial^\alpha (V^\mu(x) -W^\mu(x)) =\CO(\mu^\infty\la \Re x\ra^{-\widetilde\nu -|\alpha|}),
 $$
 uniformly with respect to  $x\in\Sigma (\mu)$ and $\mu>0$ small enough.
 \end{itemize}
 \end{prop}
 \begin{proof} We denote by $\widetilde V$ a smooth function on $\C^n$ verifying,
\begin{itemize}
\item $\widetilde V = V$ on $\R^n$;
\item For any $C>0$, one has,
$$
\overline\partial \widetilde V =\CO\left( (|\Im x|/\la\Re x\ra)^\infty \la \Re x\ra^{-\nu}\right),
$$
uniformly on $\{ |\Im x| \leq C\la \Re x\ra\}$;
\item For any $C>0$ and $\alpha\in\N^n$, one has,
$$
\partial^\alpha \widetilde V =\CO\left(  \la \Re x\ra^{-\nu-|\alpha|}\right),
$$
uniformly on $\{ |\Im x| \leq C\la \Re x\ra\}$.
\end{itemize}
Note that such a function $\widetilde V$ (called an ``almost-analytic'' extension of $V$: See, e.g., \cite{MeSj}) can easily be obtained by taking a resummation of the formal series,
\be
\label{formserie}
\widetilde V(x)\sim \sum_{\alpha\in\N^n} \frac{i^{|\alpha|}(\Im x)^\alpha}{\alpha !}\partial^\alpha V(\Re x).
\ee
Indeed, since we have $\partial^\alpha V(\Re x)=\CO(\la\Re x\ra^{-\nu-|\alpha|})$, the resummation is well defined up to $\CO\left( (|\Im x|/\la\Re x\ra)^\infty \la \Re x\ra^{-\nu}\right)$, and the standard procedure of resummation (see, e.g., \cite{DiSj, Ma1}) also gives the required estimates on the derivatives of $\widetilde V$. Conversely, by a Taylor expansion, we see that any $\widetilde V$ verifying the required conditions is necessarily a resummation of the series (\ref{formserie}).

Now, if $V^\mu$ is a $|x|$-analytic $(\mu, \widetilde\nu)$-approximation of $V$, then, for any $x=r\omega \in \Sigma (\mu)$ ($\omega\in S^{n-1}$) and $N\geq 0$, we have,
\begin{eqnarray*}
V^\mu (x) -\widetilde V (x) &=& \sum_{k=0}^{N}\frac{i^k(\Im r)^k}{k !}\partial_r^k V^\mu (\Re r\cdot\omega)\\
&&+\frac{(i\Im r)^{N+1}}{(N+1)!}\int_0^1\partial_r^{N+1}(V^\mu( (\Re r + it\Im r)\cdot\omega))dt -\widetilde V(x)\\
&=& \sum_{k=0}^{N}\frac{i^k(\Im r)^k}{k !}\partial_r^k (V^\mu (\Re x)-V(\Re x)) 
+{\cal O}(\mu^{N+1}\la \Re x\ra^{-\widetilde \nu})\\
&=& {\cal O}(\mu^{\infty}\la \Re x\ra^{-\widetilde \nu})+\CO(\mu^{N+1}\la \Re x\ra^{-\widetilde \nu}),
\end{eqnarray*}
and, similarly, for any $\alpha\in \N^n$,
$$
\partial^\alpha (V^\mu (x) -\widetilde V (x)) = \CO(\mu^{\infty}\la \Re x\ra^{-\widetilde \nu -|\alpha|}).
$$
In particular, we have proved (ii).

 Now, we proceed with the construction of such a $V^\mu$.
 
For $x\in\R^n\backslash 0$, we set $\omega = \frac{x}{|x|}$, $r=|x|$, and $s=\ln r$. In particular, for any $t$ real small enough, the dilation $x\mapsto e^t x$ becomes $(s,\omega)\mapsto (s+t,\omega)$ in the new coordinates $(s,\omega)$.

For $\omega\in S^{n-1}$ and $s\in\C$ with $|\Im s|$ small enough, we set
$\widetilde V_1(s, \omega) := \widetilde V(e^s\omega)$, where $\widetilde V$ is an almost-analytic extention of $V$ as before. Then, for $|\Im s| <2 \mu$ and $\Re s\geq -\mu$, we define,
\be
\label{44a}
 V^\mu_1(s,\omega):= \frac{e^{-\widetilde\nu s}}{2i\pi}\int_{\gamma}\frac{e^{\widetilde\nu s'}\widetilde V_1(s',\omega )}{s-s'}ds',
\ee
where $\gamma$ is the oriented complex contour,
\be
\gamma:=((+\infty, -2\mu]+2i\mu)\cup (-2\mu+2i[\mu, -\mu])\cup ([-2\mu, +\infty)-2i\mu).
\ee
Observe that, by construction, we have $\widetilde V_1(s', \omega)=\CO(e^{-\nu\Re s'})$, so that the previous integral is indeed absolutely convergent. 
Therefore, the $(s,\omega)$-smoothness and $s$-holomorphy of $ V^\mu_1$ are obvious consequences of Lebesgue's dominated convergence theorem. Since $\gamma$ is symmetric with respect to $\R$, we also have that $V^\mu_1(s,\omega)$ is real for $s$ real. Moreover, 
since  $|s-s'|\geq \mu$ on $\gamma$, we see that,
$$
 V^\mu_1(s,\omega)= \frac{e^{-\widetilde\nu s}}{2i\pi}\int_{\gamma(s)}\frac{e^{\widetilde\nu s'}\widetilde V_1(s',\omega )}{s-s'}ds' +\CO(e^{-(\nu -\widetilde\nu) /(2\mu)-\widetilde\nu\Re s}),
$$
where,
$$
\gamma(s):=\left(\gamma\cap \{ \Re s'\leq \Re s + \frac1{\mu}\} \right)\cup (\Re s + \frac1{\mu}+2i[-\mu,\mu]).
$$
In particular, $\gamma (s)$ is a simple oriented loop around $s$, and therefore, one obtains,
\begin{eqnarray}
 V^\mu_1(s,\omega)-\widetilde V_1(s,\omega)&=& \frac{e^{-\widetilde\nu s}}{2i\pi}\int_{\gamma(s)}\frac{e^{\widetilde\nu s'}\widetilde V_1(s',\omega )-e^{\widetilde\nu s}\widetilde V_1(s,\omega)}{s-s'}ds' \nonumber\\
 \label{diffV}
 &&+\CO(e^{-(\nu -\widetilde\nu) /(2\mu)-\widetilde\nu\Re s}).
\end{eqnarray}
Then,  writing,
\be
e^{\widetilde\nu s'}\widetilde V_1(s',\omega )-e^{\widetilde\nu s}\widetilde V_1(s,\omega) = (s-s')f(s,s',\omega) + \overline{(s-s')}g(s,s',\omega),
\ee
with $|\overline\partial_{s'}f | + |g|=\CO(\mu^\infty)$, 
by Stokes' formula, we see that, for $\Re s \leq 2/\mu$ and $|\Im s|\leq \mu$, we have,
$$
V^\mu_1(s,\omega)-\widetilde V_1(s,\omega) =\CO(\mu^\infty e^{-\widetilde\nu\Re s}).
$$

When $\Re s >2/\mu$  and $|\Im s|\leq \mu$,
setting,
$$
\gamma_1(s):=\left(\gamma\cap \{ \Re s'\leq  \frac1{\mu}\} \right)\cup (\frac1{\mu}+2i[-\mu,\mu]),
$$
Stokes' formula directly gives,
$$
\int_{\gamma_1(s)}\frac{e^{\widetilde\nu s'}\widetilde V_1(s',\omega )}{s-s'}ds' =\CO(\mu^\infty),
$$
and thus, using again that  $\widetilde V_1(s', \omega)=\CO(e^{-\nu\Re s'})$, in that case we obtain,
$$
| V^\mu_1(s,\omega)| + |\widetilde V_1(s,\omega)|=\CO(\mu^\infty e^{-\widetilde\nu \Re s}).
$$
In particular, in both cases we obtain,
\be
\label{estdif1}
V^\mu_1(s,\omega)-\widetilde V_1(s,\omega) = \CO(\mu^\infty e^{-\widetilde\nu\Re s}),
\ee
uniformly for $\Re s\geq -\mu$, $|\Im s|\leq \mu$ and $\mu >0$ small enough.

Then, for $\alpha\in\N^n$ arbitrary,  by differentiating (\ref{44a}) and observing that,  
\begin{eqnarray*}
&&e^{\widetilde\nu s'}\widetilde V_1(s',\omega )-\sum_{k=0}^N\frac1{k!}(s'-s)^k\partial_s^k\left( e^{\widetilde\nu s}\widetilde V_1(s,\omega)\right) \\
&& \hskip 3cm= (s'-s)^{N+1}f_N(s,s',\omega) + g_N(s,s',\omega),
\end{eqnarray*}
with $|\overline\partial_{s'}f_N|+|g_N|= \CO(\mu^\infty)$, the same procedure gives,
\be
\label{estdif2}
\partial^{\alpha}(V^\mu_1(s,\omega)-\widetilde V_1(s,\omega)) = \CO(\mu^\infty e^{-\widetilde\nu \Re s}),
\ee
uniformly for $\Re s\geq -\mu$, $|\Im s|\leq \mu$ and $\mu >0$ small enough. In particular, using the properties of $\widetilde V_1$, on the same set we also obtain,
\be
\label{estdif3}
\partial^{\alpha}V^\mu_1(s,\omega) = \CO(e^{-\widetilde\nu \Re s}),
\ee
uniformly.

Now, let $\chi_1\in C^\infty (\R ; [0,1])$ be such that $\chi_1 =1$ on $(-\infty, -1]$, and $\chi_1=0$ on $\R_+$. We set,
\be
V^\mu_2(s,\omega):= \chi_1 (s/\mu)\widetilde V_1(s,\omega) + (1- \chi_1 (s/\mu)) V^\mu_1(s,\omega).
\ee
In particular, $ V^\mu_2$ is well defined and smooth on $\R_-\cup (\R_+ +i[-\mu,\mu])$, and one has,
\begin{eqnarray*}
 &&  V^\mu_2 = \widetilde V_1 \mbox{ for  } s\in (-\infty, -\mu];\\
 &&  V^\mu_2 =  V^\mu_1 \mbox{ for  } s \in \R_+ +i[-\mu,\mu];\\
 && \partial^\alpha ( V^\mu_2 - \widetilde V_1) =  \CO(\mu^\infty ) \mbox{ for  } s\in [-\mu, \mu].
\end{eqnarray*}
Finally, setting,
\be
V^\mu (x):=  V^\mu_2( \ln |x|, \frac{x}{|x|}),
\ee
for $x\not= 0$, and $V^\mu (0) = \widetilde V(0)$,
we easily deduce from the previous discussion (in particular (\ref{estdif1}), (\ref{estdif2}) and (\ref{estdif3}), and the fact that $\partial_{r}=r^{-1}\partial_{s}$), that $V^\mu$ is a $|x|$-analytic $(\mu, \widetilde\nu)$-approximation of $V$.
  \end{proof}

\section{The analytic distortion}\label{sec-dist}
In this section, for any $\theta>0$ small enough, we construct a suitable distortion $x\mapsto x+i\theta A(x)$ verifying $A(x) =x$ for $|x|$ large enough, and such that, for $\mu\geq \theta$,  the resolvent $(P^\mu_\theta -z)^{-1}$ of the corresponding distorted Hamiltonian $P^\mu_\theta$, admits sufficiently good estimates when  $\Im z\geq h^{n_1}\theta$.

We fix $R_0\geq 1$ arbitrarily, and we have,
\begin{lem}\sl
\label{flambda}
For any $\lambda>1$ large enough, there exists $f_\lambda\in C^\infty (\R_+)$, such that,
\begin{itemize}
\item[(i)] $\supp f_\lambda \subset [R_0, +\infty)$;
\item[(ii)] $f_\lambda (r)= \lambda r$ for $r\geq 2\ln\lambda$;
\item[(iii)] $0\leq f_\lambda (r)\leq rf_\lambda'(r)\leq 2\lambda r$ everywhere;
\item[(iv)] $f_\lambda' +|f_\lambda''| = \CO(1+f_\lambda)$ uniformly;
\item[(v)] For any $k\geq 1$, $f_\lambda^{(k)}= \CO(\lambda)$ uniformly.
\end{itemize}
\end{lem}
The construction of such an $f_\lambda$ is made in Appendix \ref{app1}. 
\vskip 0.2cm
Now, we take $\lambda := h^{-n_1}$, and  we set,
\be
b (r):= \frac{1}{\lambda}f_{\lambda}(r).
\ee
By the previous lemma, $b$ verifies,
\begin{itemize}
\item $\supp b \subset [R_0, +\infty)$;
\item $b (r)=  r$ for $r\geq 2n_1\ln\frac1{h}$;
\item $0\leq b (r)\leq rb'(r)\leq 2 r$ everywhere;
\item $b' + |b''| = \CO(h^{n_1}+b)$ uniformly;
\item For any $k\geq 1$, $b^{(k)}= \CO(1)$ uniformly.
\end{itemize}
We set,
$$
A(x) := b(|x|)\frac{x}{|x|} =a(|x|)x,
$$
where $a(r):= r^{-1}b(r)\in C^\infty (\R_+)$. For $\mu\geq\theta$ (both small enough), we can define the distorted operator $P^\mu_\theta$  as in (\ref{distP}) obtained from $P^\mu$ by using the distortion,
\be
\label{defdist}
\Phi_\theta \, :\, \R^n \ni x\mapsto x+i\theta A(x) \in \C^n.
\ee
Here we use the fact that  $|A (x)|\leq 2|x|$, and we also observe that, for any $\alpha\in \N^n$ with $|\alpha|\geq 1$, one has $\partial^\alpha \Phi_\theta (x) = \CO(\theta\langle x\rangle^{1-|\alpha|})$ uniformly. 
\begin{prop} \sl
\label{estresC+}
If $R_0$ is fixed sufficiently large, then, for $0<\theta \leq\mu$ both small enough, $h>0$ small enough, $u\in H^2(\R^n)$, and $z\in\C$ such that $\Re z\in[ \lambda_0/2, 2\lambda_0]$ and $\Im z\geq h^{n_1}\theta$, one has,
$$
|\la (P^\mu_\theta -z) u,u\ra_{L^2}|\geq \frac{\Im z}2\Vert u\Vert^2_{L^2}.
$$
\end{prop}
\begin{proof} Setting $F :={}^tdA (x) = dA (x)$, and $V^\mu_\theta(x):= V^\mu (x+i\theta A (x))$, we have,
\begin{eqnarray*}
\la P^\mu_\theta u,u\ra &=& \la [(I+i\theta F (x))^{-1}hD_x]^2u,u\ra + \la V^\mu_\theta u,u\ra\\
&=& \la (1+i\theta F(x))^{-2}hD_xu, hD_xu\ra \\
&&  + ih\la [({}^t \nabla_x)(I+i \theta F(x))^{-1}](I+i \theta F(x))^{-1}h\nabla_xu, u\ra + \la V^\mu_\theta u,u\ra.
\end{eqnarray*}
Therefore, using Lemma \ref{lemA1}, and the fact that, for complex $x$, we have,
 $$
 |\Im V^\mu(x)| = \CO(|\Im x|\la \Re x\ra^{-\nu -1}),
 $$
we find,
\begin{eqnarray*}
\Im \la P^\mu_\theta u,u\ra &\leq& -\theta \Vert \sqrt{a (|x|)}\; hD_xu\Vert^2 \\
&&  + C h\theta\int \left(|b''|+\frac{b'}{|x|} + \frac{b}{|x|^2}\right) |hD_xu|\cdot |u|dx  + C_0\theta \left\Vert  \frac{\sqrt{b}}{| x|^{\frac{\nu +1}2}}u\right\Vert^2
\end{eqnarray*}
for some constants $C, C_0>0$, $C_0$ independent of the choice of $R_0$.

Thus, using the properties of $b$ after Lemma \ref{flambda}, we obtain (with some other constant $C>0$),
\begin{eqnarray}
\label{estimP1}
\Im \la P^\mu_\theta u,u\ra &\leq& -\theta\Vert \sqrt{a (|x|)}\hskip 1pt  hD_xu\Vert^2 \\
&&  + Ch\theta\int \left(|b''| + \frac{b}{\vert x\vert} + h^{n_1}\right) |hD_xu|\cdot |u|dx  + C_0R_0^{-\nu}\theta\left\Vert   \sqrt{a} u\right\Vert^2.\nonumber
\end{eqnarray}
On the other hand, for $z\in\C$, a similar computation gives,
\begin{eqnarray*}
\Re\la \sqrt{a}(P_\theta^\mu -z) u,\sqrt{a} u\ra&=&-(\Re z)\Vert  \sqrt{a}u\Vert^2\\
&& + \Re\la \sqrt{a}[(I+i \theta F (x))^{-1}hD_x]^2u,\sqrt{a}u\ra\\
&& + \Re\la \sqrt{a}V^\mu_\theta u,\sqrt{a}u\ra,\\
&\leq& -(\Re z)\Vert  \sqrt{a}u\Vert^2+ (1-2\theta)^{-2}\Vert \sqrt{a}\hskip 1pt hD_xu\Vert^2\\
&& + C h\int \left(|b''| + \frac{b}{\vert x\vert}  + h^{n_1}\right) |hD_xu|\cdot |u|dx \\
&& +C_0R_0^{-\nu}\left\Vert   \sqrt{a} u\right\Vert^2,
\end{eqnarray*}
still with $C, C_0$ positive constants, and $C_0$ independent of the choice of $R_0$. Therefore, if $\Re z\geq\lambda_0/2 >0$ and $R_0$ is chosen sufficiently large, then, for $\theta$ small enough, we obtain,
\begin{eqnarray}
\Vert  \sqrt{a}u\Vert^2&\leq& 4\lambda_0^{-1}\Vert \sqrt{a}\hskip 1pt hD_xu\Vert^2\nonumber\\
\label{estradau}
&&+ 4C\lambda_0^{-1} h\int \left(|b''| + \frac{b}{\vert x\vert}  + h^{n_1}\right) |hD_xu|\cdot |u|dx\\
&& + 4\lambda_0^{-1}|\la \sqrt{a}(P_\theta^\mu -z) u,\sqrt{a} u\ra|.\nonumber
\end{eqnarray}
The insertion of this estimate into (\ref{estimP1}) gives,
\begin{eqnarray}
\Im \la P^\mu_\theta u,u\ra &\leq& -(1-4C_0\lambda_0^{-1}R_0^{-\nu})\theta\Vert \sqrt{a }\hskip 1pt  hD_xu\Vert^2 \nonumber\\
\label{estimP2}
&&  + C'h\theta\int \left(|b''| + \frac{b}{\vert x\vert}  + h^{n_1}\right) |hD_xu|\cdot |u|dx \\
&&+ C' \theta|\la \sqrt{a}(P_\theta^\mu -z) u,\sqrt{a} u\ra|, \nonumber
\end{eqnarray}
with $C'>0$  a constant.

Now, for $r\geq 2n_1\ln\frac1{h}$, by construction we have $b'' (r)=0$, while, for $r\leq 2n_1\ln\frac1{h}$, we have,
\be
\label{estb'}
|b'' (r)| = \CO(h^{n_1} + b) = \CO(h^{n_1} + (\ln\frac1{h})a).
\ee
Then, we deduce from (\ref{estimP2}),
\begin{eqnarray}
\Im \la P^\mu_\theta u,u\ra &\leq& -(1-4C_0\lambda_0^{-1}R_0^{-\nu})\theta\Vert \sqrt{a}\hskip 1pt  hD_xu\Vert^2 \nonumber\\
\label{estimP3}
&&  + C'h\theta\ln\frac1{h}\Vert \sqrt{a }\hskip 1pt  hD_xu\Vert\cdot \Vert \sqrt{a }\hskip 1pt  u\Vert\\
&&+C'h^{n_1+1}\theta \Vert hD_xu\Vert\cdot \Vert u\Vert+ C'\theta|\la \sqrt{a}(P_\theta^\mu -z) u,\sqrt{a} u\ra|, \nonumber
\end{eqnarray}
with some other constant $C'>0$. Using again (\ref{estb'}), we also deduce from  (\ref{estradau}),
\begin{eqnarray*}
\Vert \sqrt{a }\hskip 1pt  u\Vert^2 &=& \CO( \Vert \sqrt{a }\hskip 1pt  hD_xu\Vert^2 +|\la \sqrt{a}(P_\theta^\mu -z) u,\sqrt{a} u\ra|\\
&& \hskip 4cm+ h^{n_1+1} \Vert hD_xu\Vert\cdot \Vert u\Vert),
\end{eqnarray*}
uniformly for $h>0$ small  enough,
and thus,  by (\ref{estimP3}),
\begin{eqnarray}
\label{estimP4}
\Im \la P^\mu_\theta u,u\ra &\leq& -(1-4C_0\lambda_0^{-1}R_0^{-\nu}-Ch\ln\frac1{h})\theta\Vert \sqrt{a}\hskip 1pt  hD_xu\Vert^2\\
&&+Ch^{n_1+1}\theta \Vert hD_xu\Vert\cdot \Vert u\Vert+ C\theta|\la \sqrt{a}(P_\theta^\mu -z) u,\sqrt{a} u\ra|. \nonumber
\end{eqnarray}
Finally, for $\Re z\leq 2\lambda_0$, we use the (standard) ellipticity of the second-order partial differential operator $\Re P^\mu_\theta$, and the properties of $V^\mu$, to obtain,
$$
\Re \la (P^\mu_\theta -z)u,u\ra \geq \frac1C\Vert hD_x u\Vert^2 - C\Vert u\Vert^2,
$$
where $C$ is again a new positive constant, independent of $\mu$ and $\theta$. Combining with (\ref{estimP4}), and possibly increasing the value of $R_0$, this leads to,
\begin{eqnarray}
\label{estimP5}
\Im \la (P^\mu_\theta -z)u,u\ra &\leq& (Ch^{n_1+1}\theta  -\Im z)\Vert u\Vert^2\\
&&+Ch^{n_1+1}\theta  |\la (P^\mu_\theta -z)u,u\ra|^{\frac12}\Vert u\Vert + C\theta |\la (P_\theta^\mu -z) u,u\ra|, \nonumber
\end{eqnarray}
and thus, for $\Im z\geq h^{n_1}\theta$, and  for $h,\theta >0$ small enough, we can deduce,
\be
\label{estimP6}
|\la (P^\mu_\theta -z)u,u\ra | \geq  \frac{3\Im z}4\Vert u\Vert^2-Ch^{n_1+1}\theta  |\la (P^\mu_\theta -z)u,u\ra|^{\frac12}\Vert u\Vert. 
\ee
Then, the result easily follows by solving this second-order inequation where the unkonwn variable is $ |\la (P^\mu_\theta -z)u,u\ra|^{\frac12}$, and by using again that $\Im z \gg h^{n_1+1}\theta$.
\end{proof}

%%%%%%%%%%%%%%%%SECTION 6%%%%%%%%%%%%%%%%%%%%%%%%%%%
\section{Proof of Theorem 2.1}\label{sec-th21}

%%%%%%%%%%SUBSECTION 6.2%%%%%%%%%%%%%%%%%%%%%%
 \subsection{The invertible reference operator}\label{sec-filled}

The purpose of this section is to introduce an operator without eigenvalues near $\lambda_0$, obtained as a finite-rank perturbation of $P^{\mu}_{\theta}$, $0<\theta\leq\mu$, and for which we have a nice estimate for the resolvent in the lower half plane. This operator will be used in the next section to construct a convenient Grushin problem.
 
 Let $\chi_0\in C_0^\infty (\R_+ ; [0,1])$ be equal to $1$ on $[0,1+2\lambda_0+\sup|V|]$, and let $C_0>\sup |\nabla V|$. We set,
\begin{eqnarray*}
&&  R=R(h):= 2n_1\ln\frac1{h}; \\
 && \widetilde P^{\mu}_{\theta}:= P^{\mu}_{\theta} -iC_0\theta\chi_0(h^2D_x^2 +R^{-2}x^2).
 \end{eqnarray*}
 Observe that $h^2D_x^2 +R^{-2}x^2$ is unitarily equivalent to $hR^{-1}(D_x^2+x^2)$, and therefore
 the rank of $\chi_0(h^2D_x^2 +R^{-2}x^2)$ is $\CO(R^nh^{-n})$.
 
 For $m\in\R$, we denote by $S(\la \xi\ra^m)$ the set of functions $a\in C^\infty (\R^{2n})$ such that, for all $\alpha\in\N^{2n}$, one has,
 $$
\partial^\alpha_{x,\xi}a(x,\xi ) =\CO (\la \xi\ra^m) \mbox{ uniformly}.
$$

We also denote 
\be\label{wq}
\Op(a) u (x ) = \frac1{(2\pi h)^n}\iint e^{i(x-y)\xi /h}a\big(\frac{x+y}2,\xi )u(y)dyd\xi,
\ee 
 the semiclassical Weyl quantization of such a symbol $a$.
 
 Denoting by $\widetilde p^{\mu}_{\theta}\in S(\la\xi\ra^2)$ the Weyl symbol of $\widetilde P^{\mu}_{\theta}$, we see that,
\begin{eqnarray}
\widetilde p^{\mu}_{\theta}(x,\xi )= [({}^td\Phi_{\theta}(x))^{-1}\xi]^2 + V^{\mu} (\Phi_{\theta}(x))-iC_0\theta\chi_0(\xi^2 + R^{-2}x^2) 
 \label{symbptheta}
 + \CO(h\theta\la\xi\ra ),
\end{eqnarray}
uniformly with respect to $(x,\xi)$, $\mu$, $\theta$, and $h$, and where the estimate on the remainder is in the sense of symbols (that is, one has the same estimate for all the derivatives).
In particular, we have,
\be
\label{estrep1}
 \Re\widetilde p^{\mu}_{\theta}(x,\xi) =\xi^2 +V(x)+\CO(\theta\la\xi\ra^2).
 \ee
Moreover
  \begin{itemize}
 \item If $|x|\geq R$ and $|\xi|^2\geq \lambda_0 /2$, then,
\be
\label{estimp1}
 \Im\widetilde p^{\mu}_{\theta}(x,\xi)  \leq -\frac{\theta}{C}\la\xi\ra^2 +\CO( \theta R^{-\nu}) \leq -\frac{\theta}{2C}\la\xi\ra^2;
\ee
 \item If $|x|\leq R$ and $|\xi|^2\leq 2\lambda_0+\sup|V|$, then,
\be
\label{estimp2}
 \Im\widetilde p^{\mu}_{\theta} \leq -C_0\theta + \theta\sup|\nabla V| +\CO(h\theta) \leq -\frac{\theta}{2C},
\ee
 \end{itemize}
 where $C>0$ is a constant, and the estimates are valid for $h$  small enough. (For (\ref{estimp2}), we have used the fact that $\Im [({}^td\Phi_{\theta}(x))^{-1}\xi]^2\leq 0$, that is due to the particular form of $\Phi_{\theta}(x)$. See Lemma \ref{lemA1} in appendix.)
\vskip 0.2cm
We have,
\begin{prop} \sl
\label{estresfilled}
There exists a constant $\widetilde C\geq 1$ such that, for all $\mu>0$, for all $\theta\in (0, \mu]$, for all $z$ verifying $|\Re z -\lambda_0| + \theta^{-1}|\Im z|\leq \frac4{\widetilde C}$,   and for all $h\in(0,1/\widetilde C]$, one has,
$$
\Vert (z-\widetilde P^{\mu}_{\theta})^{-1}\Vert \leq \widetilde C\theta^{-1}.
$$
\end{prop}
\begin{proof}
We take two functions $\varphi_1,\varphi_2\in C_{b}^\infty (\R^{2n} ; [0,1])$ (the space of smooth functions bounded with all their derivatives), such that,
\begin{itemize}
\item $\varphi_1^2 +\varphi_2^2 =1$ on $\R^{2n}$;
\item $\supp \varphi_1$ is included in a small enough neighborhood of $\{ \xi^2 + V(x) =\lambda_0\}$;
\item $\varphi_1 =1$ near $\{ \xi^2 + V(x) =\lambda_0\}$.
\end{itemize}
In particular, $\varphi_1$ can be chosen in such a way that, on  $\supp \varphi_1$, one has either $|x|\geq R$ together with  $|\xi|^2\geq \lambda_0/2$, or $|x|\leq R$ together with  $|\xi|^2\leq 2\lambda_0 +\sup|V|$. Therefore, we deduce from (\ref{estimp1})-(\ref{estimp2}),
$$
\frac1{\theta}\Im\widetilde p^{\mu}_{\theta}\leq -\frac1{2C} \mbox{ on } \supp\varphi_1,
$$
and thus,
\be
\label{estimp3}
\varphi_1^2\frac1{\theta}\Im\widetilde p^{\mu}_{\theta}+\frac1{2C}\varphi_1^2\leq 0 \mbox{ on } \R^{2n}.
\ee
Moreover, it is easy to check that the function $(x,\xi)\mapsto \theta^{-1}\Im\widetilde p^{\mu}_{\theta}$ is a nice symbol in $S(\la \xi\ra^2)$, uniformly with respect to $\mu$ and $\theta$, that is, for all $\alpha\in\N^{2n}$, one has,
$$
\partial^\alpha_{x,\xi}(\theta^{-1}\Im\widetilde p^{\mu}_{\theta})(x,\xi ) =\CO (\la \xi\ra^2) \mbox{ uniformly},
$$
and we see from (\ref{symbptheta}), that,
$$
\theta^{-1}\Im\widetilde p^{\mu}_{\theta}\leq CR^{-\nu} + Ch\la\xi\ra ,
$$
with some new uniform constant $C>0$.

Then, setting $\phi_j:={\rm Op}_h^W(\varphi_j)$,  writing $I= \phi_1^2u+\phi_2^2u +hQ$ where $Q$ is a uniformly bounded pseudodifferential operator, and using the sharp G\aa rding inequality, we obtain,
\begin{eqnarray*}
\la \theta^{-1}\Im\widetilde P^{\mu}_{\theta} u,u\ra &=& \la \phi_1\theta^{-1}\Im\widetilde P^{\mu}_{\theta}\phi_1u,u\ra +\la \theta^{-1}\Im\widetilde P^{\mu}_{\theta}\phi_2u,\phi_2u\ra +\CO(h\Vert u\Vert_{H^1}^2)\\
&\leq& -\frac1{2C}\Vert\phi_1u\Vert^2+CR^{-\nu}\Vert\phi_2u\Vert^2 +Ch\Vert \la hD_x\ra u\Vert^2,
\end{eqnarray*}
for all $u\in H^2(\R^n)$, and still for some new uniform constant $C>0$. Hence,
\be
\label{estimP}
 \vert \Im\la\widetilde P^{\mu}_{\theta} u,u\ra\vert \geq \frac{\theta}{2C}\Vert\phi_1u\Vert^2-C\theta R^{-\nu}\Vert\phi_2u\Vert^2 -Ch\theta\Vert \la hD_x\ra u\Vert^2.
\ee
On the other hand, since $\Re\widetilde p^{\mu}_{\theta} -\lambda_0\in S(\la \xi\ra^2)$ is uniformly elliptic on $\supp\varphi_2$, the symbolic calculus permits us to construct  $a\in S(\la\xi\ra^{-2})$ (still depending on $\mu$ and $\theta$, but with uniform estimates), such that,
$$
a\sharp (\widetilde p_{k,\theta} -\lambda_0) =\varphi_2\sharp\varphi_2 +\CO(h^\infty) \mbox{ in } S(1),
$$
where $\sharp$ stands for the Weyl composition of symbols. As a consequence, denoting by $A$ the Weyl quantization of $a$, we obtain,
$$
\Vert \la hD_x\ra\phi_2u\Vert^2= \la \la hD_x\ra^2 A (\widetilde P^{\mu}_{\theta}-\lambda_0)u,u\ra+\CO(h)\Vert u\Vert^2,
$$
and thus, 
\be
\label{estreP}
\Vert (\widetilde P^{\mu}_{\theta} -\lambda_0)u\Vert\cdot\Vert u\Vert \geq \frac1{C}\Vert \la hD_x\ra\phi_2u\Vert^2-Ch\Vert u\Vert^2.
\ee
Now, if $z\in\C$ is such that $|\Re z -\lambda_0|\leq\varepsilon$ and $|\Im z|\leq \varepsilon\theta$ ($\varepsilon >0$ fixed), we deduce from (\ref{estimP})-(\ref{estreP}),
\begin{eqnarray*}
&& \Vert(\widetilde P^{\mu}_{\theta}-z) u\Vert \cdot\Vert u\Vert
\geq 
\vert \Im\la(\widetilde P^{\mu}_{\theta} -z)u,u\ra\vert 
 \geq \frac{\theta}{2C}\Vert\phi_1u\Vert^2-C\theta R^{-\nu}\Vert\phi_2u\Vert^2 \\
&&\hskip 5cm -Ch\theta\Vert \la hD_x\ra u\Vert^2-\varepsilon \theta \Vert u\Vert^2;\\
&& \theta \Vert (\widetilde P^{\mu}_{\theta} -z)u\Vert\cdot\Vert u\Vert \geq \frac{\theta}{C}\Vert \la hD_x\ra\phi_2u\Vert^2-Ch\theta\Vert u\Vert^2-2\varepsilon\theta\Vert u\Vert^2,
\end{eqnarray*}
that yields
\begin{eqnarray}
\label{estres1}
(1+\theta)\Vert (\widetilde P^{\mu}_{\theta} -z)u\Vert\cdot\Vert u\Vert &\geq& \frac{\theta}{2C}\left(\Vert\phi_1u\Vert^2+\Vert \la hD_x\ra\phi_2u\Vert^2\right)\\
&& -\theta(2Ch +CR^{-\nu}+3\varepsilon)\Vert \la hD_x\ra u\Vert^2.\nonumber
\end{eqnarray}
Moreover, since $\xi$ remains bounded on $\supp\varphi_1$, we see that the norms $\Vert \la hD_x\ra u\Vert$ and $\Vert\phi_1u\Vert+\Vert \la hD_x\ra\phi_2u\Vert$ are uniformly equivalent, and thus, for $\varepsilon$ and $h$ small enough,   we deduce from (\ref{estres1}),
$$
\Vert (\widetilde P^{\mu}_{\theta} -z)u\Vert\cdot\Vert u\Vert\geq \frac{\theta}{4C}\Vert \la hD_x\ra u\Vert^2,
$$
and the result follows.
\end{proof}

%%%%%SUBSECTION 6.2%%%%%%%%%%%%%%%%%%%%%%%%

\subsection{The Grushin problem}\label{sec-grushin}

In this section, we reduce the estimate on $(P^\mu_{\theta}-z)^{-1}$ in Theorem \ref{Th1}, to that of a finite matrix, by means of some convenient Grushin problem.

\bigskip
Denote by $(e_1,\dots,e_M)$ an orthonormal basis of the range of the operator,
$$
K:=C_0 \chi_0(h^2D_x^2+R^{-2}x^2).
$$
 In particular, $M=M(h)$ verifies,
\be
\label{estM}
M=\CO(R^nh^{-n}).
\ee
Let $z\in\C$, and consider the two operators,
$$
R_+\, :\, \begin{array}{ccc}
L^2(\R^n)& \rightarrow &\C^M\\
u & \mapsto & (\la u,e_j\ra)_{1\leq j\leq M},
\end{array}
$$
and,
$$
R_-(z)\, :\, \begin{array}{ccc}
\C^M & \rightarrow & L^2(\R^n)\\
u^- & \mapsto & \sum_{j=1}^Mu_j^-(\widetilde P^{\mu}_{\theta}-z)e_j.
\end{array}
$$
Then,  the  Grushin operator,
$$
{\cal G}(z):=\left(\begin{array}{cc}
P^{\mu}_{\theta}-z & R_-(z)\\
R_+ & 0
\end{array}\right),
$$
is well defined from $H^2(\R^n)\times \C^M$ to $L^2(\R^n)\times \C^M$, and for $z$ as in Proposition \ref{estresfilled}, it is easy to show that ${\cal G}(z)$ is invertible, and its inverse is given by,
$$
{\cal G}(z)^{-1}:=\left(\begin{array}{cc}
E(z) & E^+(z)\\
E^-(z) & E^{-+}(z)
\end{array}\right),
$$
where,
\begin{eqnarray*}
&& E(z) = (1-T_M)(\widetilde P^{\mu}_{\theta} -z)^{-1}, \mbox{ with } T_Mv:= \sum_{j=1}^M\la v,e_j\ra e_j \,\, (v\in L^2);\\
&& E^+(z)v^+ = \sum_{j=1}^M v_j^+(e_j +i\theta(1-T_M)(\widetilde P^{\mu}_{\theta}-z)^{-1}Ke_j), \\
&& \hskip 7.5cm (v_+=(v_j^+)_{1\leq j\leq M}\in\C^M);\\
&& E^-(z)v= (\la (\widetilde P^{\mu}_{\theta}-z)^{-1}v,e_j\ra)_{1\leq j\leq M};\\
&& E^{-+}(z)v^+ = -v^+ +i\theta \left( \sum_{\ell=1}^M v_\ell^+\la (\widetilde P^{\mu}_{\theta}-z)^{-1}Ke_\ell,e_j\ra\right)_{1\leq j\leq M}.
\end{eqnarray*}
Proposition \ref{estresfilled} gives
\begin{eqnarray}
&& \Vert E(z)\Vert + \Vert E^-(z)\Vert =\CO(\theta^{-1});\\
\label{E-+bdd}
&& \Vert E^+(z)\Vert + \Vert E^{-+}(z)\Vert =\CO(1),
\end{eqnarray}
uniformly for $\mu >0$, $\theta\in (0,\mu]$, $h>0$ small enough, and  $|\Re z -\lambda_0| +\theta^{-1}|\Im z|$   small enough.

Hence, using the algebraic identity,
\be
\label{algid}
(P^{\mu}_{\theta}-z)^{-1} =E(z) - E^+(z)E^{-+}(z)^{-1}E^-(z),
\ee
we finally obtain,
\begin{prop} \sl
\label{estresgrus}
If $z\in \C$ is such that $|\Re z -\lambda_0|\leq \widetilde C^{-1}$ and  $|\Im z|\leq 2\widetilde C^{-1}\theta$, and $E^{-+}(z)$ is invertible, then so is $P^{\mu}_{\theta}-z$, and one has,
$$
\Vert (P^{\mu}_{\theta}-z)^{-1}\Vert =\CO(\theta^{-1}(1+\Vert E^{-+}(z)^{-1}\Vert )),
$$
uniformly with respect to $\mu >0$, $\theta\in (0,\mu]$, $h>0$ small enough, and  $z$ such that $|\Re z -\lambda_0|\leq \widetilde C^{-1}$ and  $|\Im z|\leq\widetilde C^{-1}\theta$.

\end{prop}
Therefore, we have reduced the study of $(P^{\mu}_{\theta}-z)^{-1}$ to that of the $M\times M$ matrix 
$E^{-+}(z)^{-1}$. 
%In the sequel, we also assume that $\widetilde C$ is chosen sufficiently large in order to have $\widetilde C^{-1}|\Im z|\leq |\arg z|\leq \widetilde C|\Im z|$ on  $\{|\Re z -\lambda_0|\leq \widetilde C^{-1}\}$.

%%%%%%%%%%SUBSECTION 6.3%%%%%%%%%%%%%%%%%%%%%%

\subsection{Using the Maximum Principle}\label{sec-maxprinc}

For $z\in J+i[-\theta/\widetilde C,2\theta/\widetilde C]$, we define,
$$
D(z) := \det E^{-+}(z).
$$
Then, $z\mapsto D(z)$ is holomorphic in $J+i[-\theta/\widetilde C,2\theta/\widetilde C]$. Using (\ref{algid})  and setting $N:=\#(\sigma(P^{\mu}_{\theta})\cap (J+i[-\theta/\widetilde C,2\theta/\widetilde C])$, we see that $D(z)$ can be written as,
$$
D(z) =G(z)\prod_{\ell=1}^{N} (z-\rho_{\ell}),
$$
with $G$ holomorphic in $J+i[-\theta/\widetilde C,2\theta/\widetilde C]$, $G(z)\not=0$ for all $z\in J-i[0,\widetilde C^{-1}\theta]$. 

Moreover, using (\ref{E-+bdd}) and (\ref{estM}), we obtain,
\be
\label{estDk}
|D(z)| =\prod_{\lambda\in\sigma (E^{-+}(z))}|\lambda|\leq \Vert E^{-+}(z)\Vert^M\leq C_1e^{C_1R^nh^{-n}},
\ee
for some uniform constant $C_1>0$.

\begin{lem}\sl
\label{lem-bord}
For every $\theta\in [0,\mu]$, there exists $r_\theta\in [\theta/(2\widetilde C),\theta/\widetilde C]$, such that for all $z\in J-ir_{\theta}$, and for all $\ell =1,\dots,N$, one has,
$$
|z-\rho_{\ell}| \geq \frac{\theta}{8\widetilde CN}.
$$
\end{lem}
\begin{proof} By contradiction, if it was not the case, then for all $t$ in $[-\theta/2\widetilde C ,-\theta/\widetilde C]$, there should exist $\ell$ such that,
$$
|t-\Im\rho_{\ell}| < \frac{\theta}{8\widetilde CN}.
$$
Therefore, the interval $[-\theta/2\widetilde C ,-\theta/\widetilde C]$ would be included in $\cup_{\ell=1}^{N}[\Im\rho_{\ell}-\theta/(8\widetilde CN), \Im\rho_{\ell}+\theta/(8\widetilde CN)]$, which is impossible because of their respective size.
\end{proof}

From now on, we assume $\CP(\widetilde{\mu},\mu;I,J)$ and setting,
$$
{\mathcal W}_\theta(J):=J+i[-r_{\theta},2\theta/\widetilde{C}],
$$
we deduce from Lemma \ref{lem-bord} that, for $\theta\in (0,\widetilde{\mu}]$,  $z$ on the boundary of 
${\mathcal W}_\theta(J)$, and for all $\ell=1,\dots,N$, we have,
$$
|z-\rho_{\ell}|\geq \frac1{C_2}\theta,
$$
for some constant $C_2>0$. As a consequence, using (\ref{estDk}), on this set we obtain,
$$
|G(z)|\leq \theta^{-C_3}e^{C_3R^nh^{-n}},
$$
with some other uniform constant $C_3>0$.
Then, the maximum principle tells us that this estimate remains valid in the interior of ${\mathcal W}_\theta(J)$, that is,
\begin{prop}\sl
\label{estG}
 There exists a constant $C_3>0$ such that,  for all $\widetilde{\mu}$, $\mu$, $I$ and $J$ verfying  (\ref{condmu}) -- (\ref{condint}) such that $\CP(\widetilde{\mu},\mu;I,J)$ holds,  one has,
$$
|G(z)|\leq  \theta^{-C_3}e^{C_3R^nh^{-n}},
$$
for all $\theta\in (0,\widetilde{\mu}]$ , $z\in {\mathcal W}_\theta(J)$, and $h\in (0, 1/C_{3}]$.
\end{prop}

%%%%%%%%%%SUBSECTION 6.4%%%%%%%%%%%%%%%%%%%%%%

 \subsection{Using the Harnack Inequality}\label{sec-harnack}

 Since $G(z)\not=0$ on ${\mathcal W}_\theta(J)$, we can consider the function,
 $$
 H(z):=  C_3R^nh^{-n}-C_3\ln\theta-\ln |G(z)|.
 $$
Then, $H$ is harmonic in ${\mathcal W}_\theta(J)$, and, by Proposition \ref{estG}, it is also nonnegative.

Using the algebraic formula,
$$
E^{-+}(z)^{-1} = -R_+(P^{\mu}_{\theta}-z)^{-1}R_-(z),
$$
and the fact that $(P^{\mu}_{\theta}-z)^{-1}R_-(z)u^-=\sum_{j=1}^Mu_j(I-i\theta (P^{\mu}_{\theta}-z)^{-1}K)e_j$, 
we deduce from Proposition \ref{estresC+} that, for $z\in [\lambda_0/2, 2\lambda_0] + i[\theta h^{n_{1}} , 1]$, one has,
$$
\Vert E^{-+}(z)^{-1}\Vert \leq 1+2C_0h^{-n_1}.
$$
As a consequence, for such values of $z$, we obtain,
$$
\frac1{D(z)}=\det E^{-+}(z)^{-1} \leq (1+2C_0h^{-n_1})^M,
$$
and thus,
$$
|G(z)| = |D(z)|\prod_{\ell=1}^{N} |z- \rho_{\ell}|^{-1}\geq \frac1{C_4} h^{n_1M},
$$
with some uniform constant $C_4>0$.
In particular, for any $\lambda \in\R$ such that $\lambda+i\theta h^{n_1}\in{\mathcal W}_\theta(J)$, this gives,
\be
\label{estHar1} 
H(\lambda +i\theta h^{n_{1}})\leq C_3R^nh^{-n}-C_3\ln\theta +\ln C_4-n_1M\ln h.
\ee
Now, the Harnack inequality tells us that, for any $\lambda , r$, such that,
$$
\dist (\lambda , \R\setminus J)\geq \widetilde C^{-1}\theta\quad;\quad r\in [0,\widetilde C^{-1}\theta)
$$
and for any $\alpha\in\R$, one has,
$$
H(\lambda + i h^{n_{1}}\theta+re^{i\alpha})\leq \frac{\widetilde C^{-2}\theta^2}{(\widetilde C^{-1}\theta-r)^2}H(\lambda+ ih^{n_{1}}\theta).
$$
In particular, setting
$$
\widetilde{\cal W}_\theta (J):=\left\{ z\in\C\,;\, \dist (\Re z , \R\setminus J)\geq \widetilde C^{-1}\theta\, , |\Im z|\leq (2 \widetilde C)^{-1}\theta\right\},
$$
and using (\ref{estHar1}),  we find,
$$
H(z)\leq 5C_3R^nh^{-n}-5C_3\ln\theta +5\ln C_4-5n_1M\ln h,
$$
for all $z\in\widetilde{\cal W}_\theta (J)$, that is,
$$
\ln |G(z)|\geq -4C_3R^nh^{-n}+4C_3\ln\theta-5\ln C_4+5n_1M\ln h,
$$
or, equivalently,
\be
\label{estG2}
|G(z)|\geq C_4^{-5}\theta^{4C_3}h^{5n_1 M}e^{-4C_3R^nh^{-n}}.
\ee
Finally, writing $E^{-+}(z)^{-1} = D(z)^{-1}\widetilde E^{-+}(z)$, where $\widetilde E^{-+}(z)$ stands for the transposed of the comatrix of $E^{-+}(z)$, we see that,
$$
\Vert E^{-+}(z)^{-1}\Vert \leq e^{CM}|G(z)|^{-1}\prod_{\ell=1}^{N} |z-\rho_{\ell}|^{-1},
$$
and therefore, we deduce from (\ref{estG2}) and (\ref{estM}),
$$
\Vert E^{-+}(z)^{-1}\Vert \leq \theta^{-C}h^{-CR^nh^{-n}}\prod_{\ell=1}^{N}  |z-\rho_{\ell}|^{-1},
$$
with some new uniform constant $C\geq 1$.
Thus, using Proposition \ref{estresgrus}, and the fact that $R = \CO(\vert\ln h\vert)$, we have proved,
\begin{prop} \sl
\label{estexpresP}
There exists a constant $\check C>0$ such that,  for all $\widetilde{\mu}$, $\mu$, $I$ and $J$ verifying  (\ref{condmu}) -- (\ref{condint}) such that $\CP(\widetilde{\mu},\mu;I,J)$ holds,  one has,
$$
\Vert (P^{\mu}_{\theta}-z)^{-1}\Vert \leq \theta^{-\check C}h^{-\check C\vert\ln h\vert^nh^{-n}}\prod_{\ell=1}^{N}  |z-\rho_{\ell}|^{-1},
$$
for all $\theta\in (0,\widetilde{\mu}]$ , $z\in \widetilde{\mathcal W}_\theta(J)$, and $h\in (0, 1/\check C]$.
\end{prop}

 %%%%%%%%%%SUBSECTION 6.5%%%%%%%%%%%%%%%%%%%%%%
  \subsection{Using the 3-lines theorem}\label{sec-3lines}
  
  Now, following an idea of \cite{TaZw}, we define,
  $$
  \Psi(z):= \int_{a}^{b}e^{-(z-\lambda)^2/\theta^2}d\lambda,
  $$
  where,
$$
[a,b]:=\{ \lambda\in\R\; ;\; \dist (\lambda , \R\setminus J)\geq \widetilde C^{-1}\theta +\check C^{1/2}\omega_{h}(\theta)\}.
$$

  We have,
  \begin{itemize}
  \item If $\Im z =2\theta h^{n_{1}}$, then,
  $$
  | \Psi(z)| \leq (b-a)e^{4h^{2n_{1}} }=\CO (h^{\delta})\leq 1;
  $$
  \item If $\Im z =-\theta/(2\widetilde C)$, then,
  $$
  | \Psi(z)|  \leq (b-a)e^{1/4\widetilde{C}^2}=\CO (h^{\delta})\leq 1;
  $$
  \item If $\Re z\in \big\{a-\check C^{1/2}\omega_{h}(\theta), b+\check C^{1/2}\omega_{h}(\theta)\big\}$ 
  and $\Im z\in [-\theta/(2\widetilde C), 2\theta h^{n_{1}}]$, then,  
  $$
  | \Psi(z)|\leq  (b-a)e^{1/4\widetilde{C}^2}e^{-\check C\omega_{h}(\theta)^2/\theta^2}=\CO (h^{\delta})\theta^{\check C}h^{\check C\vert\ln h\vert^nh^{-n}}
  \leq   \theta^{\check C}h^{\check C\vert\ln h\vert^nh^{-n}}.
  $$
  \end{itemize}
  Then, for $z\in  \widetilde{\mathcal W}_\theta(J)$, we consider the operator-valued function,
  $$
  Q(z):=  \Psi(z)\prod_{\ell=1}^{N}(z-\rho_{\ell}) (P^{\mu}_{\theta}-z)^{-1},
  $$
  that  is holomorphic on $ \widetilde{\mathcal W}_\theta(J)$ (this can be seen, e.g.,  from (\ref{algid})). Using, Proposition \ref{estresC+}, Proposition \ref{estexpresP}, and the previous properties of $ \Psi(z)$, we see that, $Q(z)$ verifies,
  \begin{itemize}
   \item If $\Im z =2\theta h^{n_{1}}$, then,
  $$
  \Vert Q(z)\Vert \leq \theta^{-1} h^{-n_{1}};
  $$
  \item If $\Im z =-\theta/(2\widetilde C)$, then,
  $$
  \Vert Q(z)\Vert \leq \theta^{-\check C}h^{-\check C\vert\ln h\vert^nh^{-n}};
  $$
 \item If $\Re z\in \{a-\check C^{1/2}\omega_{h}(\theta), b+\check C^{1/2}\omega_{h}(\theta)\}$ 
  and $\Im z\in [-\theta/(2\widetilde C), 2\theta h^{n_{1}}]$, then, 
   $$
   \Vert Q(z)\Vert \leq 1.
   $$
  \end{itemize}
Therefore, setting,
$$
 \check{\mathcal W}_\theta(J):= [a-\check C^{1/2}\omega_{h}(\theta), b+\check C^{1/2}\omega_{h}(\theta)] +i[-\theta/(2\widetilde C), 2\theta h^{n_{1}}],
$$
 (that is included in $\widetilde{\mathcal W}_\theta(J)$), we see that the subharmonic function $z\mapsto \ln\Vert Q(z)\Vert$ verifies,
 $$
 \ln \Vert Q(z) \Vert \leq \psi (z) \,\, \mbox{ on } \partial\check{\mathcal W}_\theta(J),
 $$
 where $\psi $ is the harmonic function defined by,
\begin{eqnarray*}
 \psi (z):&=& \frac{2\theta h^{n_{1}}-\Im z}{2\theta h^{n_{1}}+\theta/(2\widetilde C)}\;\check C(\vert\ln h\vert^{n+1}h^{-n}+\vert\ln\theta\vert) \\
 &&+\frac{\Im z+\theta/(2\widetilde C)}{2\theta h^{n_{1}}+\theta/(2\widetilde C)}\;\vert\ln(\theta h^{n_{1}})\vert.
\end{eqnarray*}
 As a consequence, by the properties of subharmonic functions, we deduce that $\ln\Vert Q(z)\Vert\leq \psi (z)$ everywhere in $\check{\mathcal W}_\theta(J)$, and in particular, for $|\Im z|\leq 2\theta h^{n_{1}}$, we obtain,
 $$
 \ln\Vert Q(z)\Vert\leq 8\widetilde C\check C h^{n_1}(\vert\ln h\vert^{n+1}h^{-n}+\vert\ln\theta\vert)+\vert \ln(\theta h^{n_{1}})\vert
 $$
Hence, since $n_1>n$, we have proved the existence of some uniform constant $C\geq 1$, such that,
$$
\ln\Vert Q(z)\Vert\leq \ln C+C\vert \ln(\theta h^{n_{1}})\vert \quad \mbox{ for } z\in\check{\mathcal W}_\theta(J) \mbox{ and } h\in(0, 1/C].
$$
Coming back to $P^{\mu}_{\theta}$, this means that, for $ z\in\check{\mathcal W}_\theta(J)$ different from $\rho_{1},\dots,\rho_{N}$, we have,
$$
|\Psi(z)|\,\Vert  (P^{\mu}_{\theta}-z)^{-1}\Vert \leq C(\theta h^{n_{1}})^{-C}\prod_{\ell=1}^{N}|z-\rho_{\ell}|^{-1}.
$$
On the other hand, if 
$\dist (\Re z , \R\setminus J)\geq 2\check C^{1/2}\omega_{h}(\theta)$, and $|\Im z|\leq 2\theta h^{n_{1}}$, then,  writing $z= s+it$, we see that,
$$
\Psi (z) = \theta e^{t^2/\theta^2}\int_{(a-s)/\theta}^{(b-s)/\theta}e^{-u^2+2i(t/\theta)u}du.
$$
Now,  $\vert t/\theta\vert\leq 2 h^{n_{1}}\rightarrow 0$ uniformly, and we see that,
$$
(a-s)/\theta \leq \widetilde{C}^{-1}-\check C^{1/2}\omega_{h}(\theta)/\theta \leq \widetilde{C}^{-1}-(h^{-n}\vert \ln h\vert )^{1/2} \rightarrow -\infty, \mbox{ uniformly}.
$$
 In the same way, we have $(b-s)/\theta\rightarrow +\infty$ uniformly as $h\rightarrow 0_+$. Therefore, we easily conclude that,
$$
|\Psi (z)|\geq \frac{\theta}{C},
$$
when $h\in (0,1/C]$, with some new uniform constant $C>0$.

\bigskip
As a consequence, using also that $\theta\leq h^{\delta}$, we finally obtain,
\begin{prop}\sl 
\label{Prop9.1}There exists a constant $C_0\geq1$, such that, for all $\widetilde{\mu}$, $\mu$, $I$ and $J$ verifying  (\ref{condmu}) -- (\ref{condint}), the property $\CP(\widetilde{\mu},\mu;I,J)$ implies,
\be
\label{resolvest}
\Vert (P^{\mu}_{\theta} -z)^{-1}\Vert \leq C_0\theta^{-C_0}\prod_{\ell=1}^{N}|z-\rho_{\ell}|^{-1},
\ee
for all $z\in J'+i[-2\theta h^{n_{1}}, 2\theta h^{n_{1}}]$, and for all $h\in (0,1/C_0]$, where,
$$
J'= \{\lambda\in\R\; ;\; \dist (\lambda , \R\setminus J)\geq C_{0}\omega_{h}(\theta)\}.
$$
\end{prop}
Since $J' =J+\CO (\omega_{h}(\theta))$, Theorem \ref{Th1} is proved.

%%%%%%%%%%SECTION 7%%%%%%%%%%%%%%%%%%%%%%

\section{Proof of Theorem 2.2}\label{sec-red}

Suppose $\CP (\widetilde{\mu},\mu; I,J)$ holds, and $\widetilde{\mu}\geq \mu^{N_{0}}$ for some constant $N_{0}\geq 1$. Then, for  any $\theta\in [\mu^{N_{0}},\widetilde{\mu}]$, any constant $N_{1}\geq 1$, and any $\mu'\in [\max (\theta,\mu^{N_{1}}), \mu^{1/N_{1}}]$,  we can write,
\be
\label{pertres}
z-P^{\mu'}_{\theta} = (z-P^{\mu}_{\theta})(1+(z-P^{\mu}_{\theta})^{-1}W),
\ee
with,
\begin{eqnarray}
W:= P^{\mu}_{\theta}-P^{\mu'}_{\theta}=V^{\mu} (x+iA_{\theta} (x))-V^{\mu'} (x+iA_{\theta} (x))
\label{estWmu}
= \CO(\mu^\infty \la x\ra^{-\nu}),
\end{eqnarray}
uniformly (see Section \ref{sec-prelim}). Moreover, taking $J'$ as in Proposition \ref{Prop9.1}, we have,
\begin{lem}\sl
\label{lem-estbord} Let $N\geq 1$ be a constant, such that $N\geq \#\Gamma (\widetilde{\mu},\mu,J)$ for all $h$ small enough.
Then, for  any $\theta\in [\mu^{N_{0}},\widetilde{\mu}]$, there exists $\tau\in [\theta h^{n_{1}}, 2\theta h^{n_{1}}]$, such that,
\be
\label{estbord}
\dist (\partial (J'+i[-\tau, \tau]), \Gamma (\widetilde{\mu},\mu,J))\geq \frac{\theta h^{n_{1}}}{4N}.
\ee
Here, $\partial(J'+i[-\tau, \tau])$ stands for the boundary of $J'+i[-\tau, \tau]$.
\end{lem}

\begin{proof} If it were not the case, using $\CP (\widetilde{\mu},\mu; I,J)$, we see that, for all $t\in  [-2\theta h^{n_{1}}, - \theta h^{n_{1}}]$, there should exist $\rho \in \Gamma (\widetilde{\mu},\mu,J)$, such that,
$$
|t -\Im \rho| \leq \frac{\theta h^{n_{1}}}{4N}.
$$
That is, we would have,
$$
[-2\theta h^{n_{1}}, - \theta h^{n_{1}}]\subset \bigcup_{\rho\in \Gamma (\widetilde{\mu},\mu,J)}\big[\rho - \frac{\theta h^{n_{1}}}{4N}, \rho+\frac{\theta h^{n_{1}}}{4N}\big],
$$
which, again, is not possible because of the respective size of these two sets.\end{proof}
\begin{rem}\sl
With a similar proof, we see that the result of Lemma \ref{lem-estbord} remains valid if one replaces the interval $[\theta h^{n_{1}}, 2\theta h^{n_{1}}]$ by $[\theta h^{n_{1}}, \theta h^{n_{1}} + (\theta h^{n_{1}})^{M}]$, and one substitutes $(\theta h^{n_{1}})^{M}$ to $\theta h^{n_{1}}$ in (\ref{estbord}), where $M\geq 1$ is any arbitrary fixed number.
\end{rem}

Using Lemma \ref{lem-estbord} and Theorem \ref{Th1}, we see that, for any $z\in\partial(J'+i[-\tau, \tau])$, we have,
$$
\Vert (P^{\mu}_{\theta} -z)^{-1}\Vert \leq C_1\theta^{-C_1}\leq C_{1}\mu^{-C_{1}N_{0}},
$$
with some new uniform constant $C_1$, and thus, by (\ref{pertres}) and (\ref{estWmu}),
 $z-P^{\mu'}_{\theta}$ is invertible, too, for $z\in\partial(J'+i[-\tau, \tau])$, and the two spectral projectors,
\begin{eqnarray}
&&  \Pi:=\frac1{2i\pi}\oint_{\partial(J'+i[-\tau, \tau])} (z-P^{\mu}_{\theta})^{-1}dz;\nonumber\\
\label{projpert}
&&{}\\
\nonumber
&& \Pi':=\frac1{2i\pi}\oint_{\partial(J'+i[-\tau, \tau])} (z-P^{\mu'}_{\theta})^{-1}dz,
\end{eqnarray}
are well-defined and verify,
\be
\label{compPiPihat}
 \Vert \Pi- \Pi'\Vert = \CO(\mu^\infty).
\ee
In particular,  $\Pi$ and $\Pi'$ have the same rank ($\leq N$), and one has,
\be
\label{compPPhat}
\Vert P^{\mu}_{\theta}\Pi-P^{\mu'}_{\theta} \Pi'\Vert = \CO(\mu^\infty).
\ee
Therefore, by standard finite dimensional arguments, the two sets $\sigma (P^{\mu'}_{\theta})\cap(J'+i[-\tau, \tau])$ and $\sigma (P^{\mu}_{\theta})\cap(J'+i[-\tau, \tau])$ coincide up to $\CO(\mu^\infty)$ uniformly, and Theorem \ref{Th2} follows.

%%%%%%%%%SECTION 8%%%%%%%%%%%%%%%%%%%

\section{Proof of Theorem 2.5}\label{proofTh3}

Now, for any integer $k\geq 0$, we set,
$$
\mu_k := h^{kn_1}\widetilde{\mu}.
$$
Since $\CP (\widetilde{\mu},\mu; I,J)$ holds, we can apply Theorem \ref{Th2} with $\mu'\in [\mu_{1}, \mu_{0}]$, and deduce the existence of $J_{1}\subset J$, with $J_{1}=J+\CO (\omega_{h}(\mu_{0}))$ and $I_{1}\supset I$ with $I_{1}=I+\CO (\mu_{0}^\infty)$, independent of $\mu'$, such that, $\CP (h^{n_{1}}\mu',\mu'; I_{1},J_{1})$ holds. In particular,  $\CP (\mu_{1},\mu_{0}; I_{1},J_{1})$ holds, and we can apply Theorem \ref{Th2} again, this time with $\mu'\in [\mu_{2}, \mu_{1}]$. Iterating the procedure, we see that, for any $k\geq 0$, there exists,
$$
 I_{k+1}= I_{k}+\CO (\mu_{k}^\infty), 
\  J_{k+1} = J_{k}+\CO (\omega_{h}(\mu_{k})) 
$$
hence,
$$
I_{k+1}= I+\CO (\mu_{0}^\infty +\cdots+\mu_{k}^\infty),\ 
J_{k+1} = J+\CO (\omega_{h}(\mu_{0})+\cdots +\omega_{h}(\mu_{k})),
$$
where the $\CO$'s are also uniform with respect to $k$, such that $\CP (h^{n_{1}}\mu',\mu'; I_{k+1},J_{k+1})$ holds for all $\mu'\in [\mu_{k+1}, \mu_{k}]$.

Since the two series $\sum_{k}\omega_{h}(\mu_{k})=\CO(\omega_{h}(\widetilde\mu))$ and $\sum_{k}\mu_{k}^M=\CO(\mu^M)$ ($M\geq 1$ arbitrary) are convergent, one can find $I'=I+\CO (\mu^\infty)$ and $J'=J+\CO (\omega_{h}(\widetilde{\mu}))$, such that,
$$
I'\supset \bigcup_{k\geq 0}I_{k}\quad ;\quad J'\subset \bigcap_{k\geq 0}J_{k}.
$$
 Then, by construction, $\CP (h^{n_{1}}\mu',\mu'; I',J')$ holds for all $\mu' \in (0, \widetilde{\mu}]$, and Theorem \ref{Th3} is proved.

%%%%%%%%%%%%%%%%SECTION 9%%%%%%%%%%%%%%%%%%%%%%%

\section{Proof of Theorem 2.6 -- The set of resonances}
\label{cons}

From the proof of Theorem \ref{Th3} (and with the same notations) we learn that, for all $k\geq 0$,  $\CP(\mu_{k+1}, \mu_{k};I_{k+1},J_{k+1})$ holds. Therefore, applying Theorem \ref{Th2} with $\theta = \mu' = \mu_{k+1}$, we obtain that there exist $\tau_{k+2}\in [\mu_{k+2},2\mu_{k+2}]$, $J_{k+1}' = J_{k+1} +\CO (\omega_{h}(\mu_{k+1}))$, and  a bijection,
$$
b_{k} : \Gamma(P^{\mu_{k}})\cap (J_{k+1}'-i[0,\tau_{k+2}])\to \Gamma(P^{\mu_{k+1}})\cap (J_{k+1}'-i[0,\tau_{k+2}])
$$
such that,
\be
\label{estb-id}
b_{k}(\lambda ) -\lambda = \CO (\mu_{k}^\infty ) \mbox{ uniformly.}
\ee
In addition, we deduce from the proof of Theorem \ref{Th2} (in particular Lemma \ref{lem-estbord}), that the $\tau_{k}$'s can be chosen in such a way, that,
\be
\label{estbordJk}
\dist (\partial (J_{k+1}'+i[-\tau_{k+2}, \tau_{k+2}]), \Gamma (P^{\mu_{k}}))\geq \frac{\mu_{k}^C}{C},
\ee
for some constant $C>0$.
Setting
$$
\Lambda_k :=  \Gamma(P^{\mu_{k}})\cap (J_{k+1}'-i[0,\tau_{k+2}]),
$$
where the elements are repeated according to their multiplicity,
and, starting from an arbitrary element $\lambda_j$ of $\Lambda_0$ ($1\leq j\leq N:=\# \Lambda_0 =\CO (1)$),  we distinguish two cases.

\bigskip
\noindent
{\bf -Case A:} For all $k\geq 0$, $b_k\circ b_{k-1}\circ \dots \circ b_0(\lambda_j)\in \Lambda_{k+1}$.

In that case, we can consider the sequence defined by,
$$
\lambda_{j,k}:=b_k\circ b_{k-1}\circ \dots \circ b_0(\lambda_j),
$$
($k\geq 0$), and, using (\ref{estb-id}), we see that, for any $k>\ell\geq 0$, we have,
$$
|\lambda_{j,k}-\lambda_{j,\ell}|\leq \sum_{m=\ell}^{k-1}|\lambda_{j,m+1}-\lambda_{j,m}|\leq C_1\sum_{m=\ell}^{k-1}\mu_{m+1}\leq C_1\mu_0\frac{h^{n_1\ell}}{1-h^{n_1}},
$$
so that $(\lambda_{j,k})_{k\geq 1}$ is a Cauchy sequence, and we set,
$$
\rho_j:= \lim_{k\rightarrow +\infty}\lambda_{j,k}.
$$
Notice that according to this definition, we have a natural notion of multiplicity of a resonance $\rho$, namely the number of sequences $\rho_{j}$ converging to $\rho$.

\bigskip
\noindent
{\bf -Case B: } There exists $k_j\geq 0$ such that $b_{k-1}\circ \dots \circ b_0(\lambda_j)\in\Lambda_k$ for all $k\leq k_j$, while $b_{k_j}\circ  \dots \circ b_0(\lambda_j)\notin \Lambda_{k_{j}+1}$.
(Here, and in the sequel, we use the convention of notation: $b_{-1} \circ b_0:=Id$.)

Then, we set,
$$
\rho_j:=  b_{k_j}\circ  \dots \circ b_0(\lambda_j).
$$

In particular, since, by construction, $\Re\rho_{j}\in I_{k_{j}+2}\subset J_{k_{j}+1}$, and $\rho_{j}\notin \Lambda_{k_{j}+1}$, we see that, necessarily,  $\Im\rho_j\in [- \tau_{k_j+2},-\tau_{k_j+3})$. 

\bigskip

Moreover, if, in Case A, we set $k_j:=+\infty$, then, for any $j=1,\dots ,\#\Lambda_{0}$ and $k\geq 0$, in both cases we have the equivalence,
\be
\label{rhoj}
\vert\Im\rho_j\vert\leq \tau_{k+2} \, \Longleftrightarrow \, k\leq k_j.
\ee

\vskip 0.2cm
Now, if $\mu' \in (0,\widetilde{\mu}]$, then $\mu' \in (\mu_{k+1}, \mu_k]$ for some $k\geq 0$, and Theorem \ref{Th2} tells us that $\Gamma (P^{\mu'})\cap (J_{k+1}'-i[0,\tau_{k+2}])$ coincides with  $\Lambda_k$ up to $\CO(\mu_{k}^\infty)$ ($= \CO((\mu')^\infty)$). Therefore, setting,
$$
\Lambda:=\{ \rho_1,\dots ,\rho_N\},
$$ 
the first part of Theorem \ref{Th4} will be proved if we can show the existence, for any $k\geq 0$, of a bijection,
$$
\widetilde b_k \, :\, \Lambda \cap(J_{k+1}'-i[0,\tau_{k+2}]) \rightarrow \Lambda_k ,
$$
such that $\widetilde b_k(\rho) -\rho = \CO(\mu_k^\infty)$ uniformly. (Actually, we do not necessarily have $\tau_{k+2}\in [h^{2n_{1}}\mu', 2h^{2n_{1}}\mu']$, but rather, $\tau_{k+2}\in [h^{2n_{1}}\mu', 2h^{n_{1}}\mu')$. However, if $\tau_{k+2}\geq 2h^{2n_{1}}\mu'$, an argument similar to that of Lemma \ref{lem-bord} or Lemma \ref{lem-estbord} gives the result stated in Theorem \ref{Th4}.)

By construction, we have,
$$
\Lambda_{k}=\{ b_{k-1}\circ\dots \circ b_0 (\lambda_j)\,; \, j=1,\cdots, N \mbox{ such that } k_j\geq k\}.
$$
while, by (\ref{rhoj}),
$$
\Lambda \cap (J_{k+1}'-i[0,\tau_{k+2}])=\{\rho_j\, ;\,  j=1,\cdots, N \mbox{ such that } k_j\geq k\}.
$$
Then, for all $j$ verifying $k_j\geq k$, we set,
$$
\widetilde b_{k} (\rho_j):= b_{k-1}\circ\dots \circ b_0 (\lambda_j),
$$
so that $\widetilde b_k$ defines a bijection between $\Lambda \cap (J_{k+1}'-i[0,\tau_{k+2}])$ and $\Lambda_{k}$. Moreover, in Case A, for any $M\geq 1$, we have,
$$
|\widetilde b_{k} (\rho_j) -\rho_j|=\lim_{\ell\rightarrow +\infty}|b_{\ell}\circ\dots\circ b_k (\widetilde b_k(\lambda_j))-\widetilde b_k(\lambda_j)|\leq \sum_{m=k}^{+\infty }C_M\mu_m^M= \frac{C_M\mu_k^M}{1-h^{n_1}},
$$
while, in Case B, we obtain,
$$
|\widetilde b_{k} (\rho_j) -\rho_j|=|b_{k_j}\circ\dots\circ b_k (\widetilde b_k(\lambda_j))-\widetilde b_k(\lambda_j)|\leq \sum_{k\leq m\leq k_j}C_M\mu_m^M\leq \frac{C_M\mu_k^M}{1-h^{n_1}},
$$
(with the usual convention $\sum_{m\in\emptyset} :=0$). Therefore, in both  cases, for $h>0$ small enough, we find,
$$
|\widetilde b_{k} (\rho_j) -\rho_j|\leq 2C_M\mu_k^M,
$$
and this gives the first part of Theorem \ref{Th4}.

\bigskip

Concerning the second part of Theorem \ref{Th4}, let $\widetilde\Lambda$ be another set verifying ($\star$). 
In particular, for any $k\geq 0$, there exist $\tau_{k+2}, \widetilde\tau_{k+2} \in [\mu_{k+2}, 2\mu_{k+2}]$, such that,  $\widetilde\Lambda\cap(J_{k+1}'-i[0,\widetilde\tau_{k+2}])$ (resp.  $\Lambda\cap(J_{k+1}'-i[0,\tau_{k+2}])$) coincides with,
$\widetilde\Lambda_k:= \Gamma (P^{\mu_k})\cap  (J_{k+1}'-i[0,\widetilde\tau_{k+2}])$,
(resp. $\Lambda_k$), up to $\CO(\mu_{k}^\infty)$.

Therefore, taking $k=0$, and using again an argument similar to that of Lemma \ref{lem-bord} or Lemma \ref{lem-estbord}, that gives the existence of $\tau'\in [\frac12\mu_{2},\mu_{2}]$ and $C>0$ constant, such that,
\be
\label{estbordtau'}
\dist (\partial (J_{1}'+i[-\tau', \tau']), \Gamma (P^{\mu_{0}}))\geq \frac{\mu_{0}^C}{C},
\ee
we obtain that
 the two sets $\Lambda\cap(J_{1}'-i[0,\tau'])$ and $\widetilde \Lambda\cap(J_{1}'-i[0,\tau'])$ coincide up to  $\CO(\mu_{0}^\infty)$, and thus have same cardinal. For $k\geq 0$, we denote by,
\begin{eqnarray*}
&& B_k \, :\, \Lambda_k\, \rightarrow \, \Lambda\cap(J_{k+1}'-i[0,\tau_{k+2}]);\\
&& \widetilde B_k \, :\, \widetilde\Lambda_k\, \rightarrow \, \widetilde\Lambda\cap(J_{k+1}'-i[0,\widetilde\tau_{k+2}])\},
\end{eqnarray*}
the corresponding  bijections. Then, thanks to (\ref{estbordtau'}), we can consider the bijection,
$$
\varphi_0= \widetilde B_0\circ B_0^{-1}{}_{\left\vert_{\Lambda\cap(J_{1}'-i[0,\tau'])}\right.} \, :\, \Lambda\cap(J_{1}'-i[0,\tau'])\rightarrow \widetilde\Lambda\cap(J_{1}'-i[0,\tau']):
$$
Using (\ref{estbordJk}) and the fact that  $\widetilde B_{k}$ differ from the identity by $\CO (\mu_{k}^\infty)$, we see that, for $k\geq1$,
\be
\label{estborddeuxJk}
\dist (\partial (J_{k+1}'+i[-\tau_{k+2}, \tau_{k+2}]), \widetilde\Lambda)
\geq \frac{\mu_{k}^C}{C},
\ee
for some other constant $C>0$.

Then setting
$$
{\cal E}_0 := \Lambda\cap \{ -\tau'\leq\Im z< -\tau_{3}\},
$$
and, for $k\geq 1$,
$$
{\cal E}_k := \Lambda\cap \{ -\tau_{k+2}\leq\Im z< -\tau_{k+3}\}.
$$
we see that, for all $k\geq 1$, the application,
\be
\label{bijtran}
\widetilde B_{k}\circ  B_{k}^{-1}{}_{\left\vert_{{\cal E}_k}\right.} : {\cal E}_k \to \widetilde\Lambda\cap \{ -\tau_{k+2}\leq\Im z< -\tau_{k+3}\},
\ee
is a bijection.

\bigskip

Finally, for $\rho \in \Lambda\cap(J_{1}'-i[0,\tau'])$, we define,
\begin{itemize}
\item $B (\rho) =\widetilde B_{k}\circ  B_{k}^{-1}(\rho)$, if  $\rho\in{\cal E}_k$ for some $k\geq 0$;
\item $B (\rho) =\rho$, if $\rho\in\R$.
\end{itemize}
We first show,
\begin{lem} 
\label{lambdareel}
$\Lambda\cap\R =\widetilde\Lambda\cap\R$.
\end{lem}
\begin{proof}
We only show that any $\rho$ in $\Lambda\cap\R$ is also in $\widetilde\Lambda$, the proof of the other inclusion being similar. For such a $\rho$,  $B_k^{-1}(\rho)\in \Lambda_k $ is well defined for all $k\geq 1$, and since $B_k^{-1}$ differs from the identity by $\CO(\mu_k^\infty)$, we obtain,
$$
\alpha_k:= B_k^{-1}(\rho) \rightarrow \rho \quad {\rm as }\quad k\rightarrow +\infty.
$$
On the other hand, since $\Lambda_{k+1}\subset \widetilde\Lambda_k= \widetilde B_k^{-1} (\widetilde\Lambda )$, there exists some $\widetilde\rho_k\in\widetilde\Lambda$ such that $\alpha_{k+1} =\widetilde B_k^{-1}(\widetilde\rho_k)$. By taking a subsequence, we can  assume that $\widetilde\rho_k$ admits a limit $\widetilde \rho \in\widetilde\Lambda$ as $k\rightarrow +\infty$. Then, using that $\widetilde B_k^{-1}$ differs from the identity by $\CO(\mu_k^\infty)$, we also obtain,
$$
\alpha_{k+1} \rightarrow \widetilde \rho \quad {\rm as }\quad k\rightarrow +\infty.
$$
Therefore, we deduce that $\rho =\widetilde\rho\in \widetilde\Lambda$ and the lemma is proved.
\end{proof}
Using Lemma \ref{lambdareel}, we see that  the application $B$ is  well defined  from $\Lambda\cap(J_{1}'-i[0,\tau'])$ to $\widetilde\Lambda\cap(J_{1}'-i[0,\tau'])$. Moreover, if  $\rho\in{\cal E}_k$  for some $k\geq 0$,  we have,
$$
|B(\rho) -\rho| = |\widetilde B_{k}\circ  B_{k}^{-1} (\rho)-\rho| = \CO(\mu_k^\infty ),
$$
and, since $\tau_{k+3}\leq |\Im \rho|\leq \tau_{k+2} = \CO(h^{2n_1})$, we also have,
$$
\mu_{k}\leq h^{-3n_1}\tau_{k+3} \leq h^{-3n_1} |\Im \rho|\leq C|\Im\rho |^{1/C},
$$
where $C>0$ is a large enough constant.
Thus, we always have, 
$$
|B(\rho) -\rho|= \CO(|\Im\rho|^\infty).
$$

Therefore,  it just remains to
see that $B$ is a bijection, but this is an obvious consequence of (\ref{bijtran}), Lemma \ref{lambdareel}, and the definition of $B$. Thus Theorem \ref{Th4} is proved.

%%%%%%%%%%%%%%%%%%%%%SECTION 10%%%%%%%%%%%%%%%%%%%

\section{Shape resonances}\label{sec-exemp}

Here we prove Theorem \ref{Th5}.
Under the assumptions of Section \ref{exemp}, one can construct, as in \cite{GeMa}, a function $G_{1}\in C^\infty (\R^{2n})$, supported near $p^{-1}([\lambda_{0}-2\varepsilon, \lambda_{0}+2\varepsilon])\backslash\{x_{0}\}$ for some $\varepsilon >0$, such that, 
\begin{align}
\label{Galinfini}
& G_{1}(x,\xi ) = x\cdot \xi \mbox { for } x \mbox{ large enough, } \vert p(x,\xi )-\lambda_{0}\vert\leq \varepsilon;\\
\label{Gfuite}
& H_{p}G_{1}(x,\xi ) \geq \varepsilon \mbox{ for } x\in \R^n\backslash \mathrm{\ddot O} \mbox{ and }  \vert p(x,\xi )-\lambda_{0}\vert\leq \varepsilon.
\end{align}
We also set,
$$
\widetilde{P}:= P+W,
$$
where $W=W(x)\in C^\infty (\R^n)$ is a non negative function, supported in a small enough neighborhood of $x_{0}$, and such that $W(x_{0})>0$. In particular, denoting by $\widetilde{p}(x,\xi )=\xi^2 +V(x)+W(x)$ the principal symbol of $\widetilde{P}$, we have $\widetilde{p}^{-1}(\lambda_{0})\subset (\R^n\backslash \mathrm{\ddot O})\times\R^n$, and thus $\lambda_{0}$ is a non-trapping energy for $\widetilde{P}$.

\bigskip

Now, we take $\mu$ and $\widetilde\mu$ such that, 
$$
\mu \leq h^\delta\quad ; \quad \widetilde{\mu}\leq \min (\mu, h^{2+\delta})
$$
with $\delta >0$  arbitrary (so that $\mu,\widetilde{\mu}$
 verify (\ref{condmu})), and we denote by $V^\mu$ a $|x|$-analytic $(\mu, \widetilde\nu)$-approximation of $V$ as before. We also set,
$$
P^\mu = -h^2\Delta + V^\mu\quad ;\quad \widetilde{P}^\mu = P^\mu +W,
$$
and, if in (\ref{Utheta}) we take $A$ supported away from $\supp W$, we see that the distorted operators $P^\mu_{\theta}$ and $\widetilde{P}^\mu_{\theta}$ are well defined for $0<\theta\leq\widetilde{\mu}$. Then, we set,
$$
G(x,\xi ):= G_{1}(x,\xi ) -A(x)\cdot\xi,
$$
that, by (\ref{Galinfini}), is in $C_{0}^\infty (\R^n ;\R)$, and we consider its semiclassical Weyl-quantization $G^W=\Op(G)$ (see \ref{wq}).

Since $\theta /h^2\leq \widetilde{\mu}/h^2\leq h^\delta$, a straightforward computation shows that the operator,
$$
R^\mu_{\theta}:= \frac1\theta\Im\left(e^{\theta G^W/h}\widetilde{P}^\mu_{\theta}e^{-\theta G^W/h}\right)
$$
is a semiclassical pseudodifferential operator, with symbol $r^\mu_{\theta}$ verifying,
\begin{align*}
& \partial^\alpha  r^\mu_{\theta} =\CO (\la\xi\ra^2) \mbox{ for all } \alpha\in\N^{2n} ;\\
& r^\mu_{\theta}(x,\xi ) = -H_{\widetilde{p}^{\widetilde{\mu}}}(A(x)\cdot\xi + G)+\CO (h^\delta)=-H_{p}G_{1}(x,\xi )+\CO (h^\delta),
\end{align*}
uniformly with respect to $\theta\in (0,\widetilde{\mu}]$ and $h>0$ small enough. As a consequence, using (\ref{Gfuite}), we see that $R^\mu_\theta$ is elliptic in a neighborhood of $\{ p(x,\xi )+W(x)=\lambda_0\}$ (uniformly with respect to $\theta$ and $\mu$). Then, by arguments similar to those of Section \ref{sec-filled}, we deduce that the operator
$$
Q_\theta^\mu := e^{\theta G^W/h}\widetilde{P}^\mu_{\theta}e^{-\theta G^W/h}
$$
verifies
$$
\Vert (Q_\theta^\mu -z)^{-1}\Vert =\CO (\theta^{-1}),
$$
uniformly for $\vert \Re z -\lambda_0\vert + \theta^{-1}\vert \Im z\vert$ small enough, $\theta\in (0,\widetilde\mu]$, and $h>0$ small enough. Since $\Vert \theta G^W /h\Vert \to 0$ uniformly as $h\to 0$, this also gives,
$$
\Vert (\widetilde P_\theta^\mu -z)^{-1}\Vert =\CO (\theta^{-1}),
$$
and from this point, one can follow all the procedure used in \cite{HeSj} Sections 9 and 10. In particular, using the same notations as in \cite{HeSj}, by Agmon-type inequalities we see that the distribution kernel $K_{(\widetilde P_\theta^\mu -z)^{-1}}$ of $(\widetilde P_\theta^\mu -z)^{-1}$ verifies, 
$$
K_{(\widetilde P_\theta^\mu -z)^{-1}}(x,y) =\widetilde\CO(\theta^{-1}e^{-d(x,y)/h})
$$
where $d(x,y)$ stands for the Agmon distance between $x$ and $y$ (see  \cite[Lemma 9.4]{HeSj}). Then, assuming $\theta =\tilde\mu \geq e^{-\eta /h}$ for some $\eta >0$ constant small enough, and performing a suitable Grushin problem as in \cite{HeSj}, we deduce that the resonances of $P^\mu$ in $[\lambda_0, \lambda_0+ Ch]-i[0, \lambda_0\min (\mu, h^{2+\delta})]$ ($C>0$ constant arbitrary) are close to the eigenvalues of the Dirichlet realization of $P$ on $\{ d(x,\R^n\backslash\mathrm{\ddot O})\geq \eta /3)\}$, up to $\CO (e^{-2(S_0-\eta)/h})$. Since these eigenvalues are real and admit semiclassical asymptotic expansions of the form,
$$
\lambda_k\sim \lambda_0 + e_kh +\sum_{\ell\geq 1}\lambda_{k,\ell}h^{1+\frac\ell{2}}
$$
(where the $e_k$'s are as in Theorem \ref{Th5}), we obtain for the corresponding resonances $\rho_k$ of $P^\mu$,
\be
\label{estresF}
\Re\rho_k\sim \lambda_0 + e_kh +\sum_{\ell\geq 1}\lambda_{k,\ell}h^{1+\frac\ell{2}}\quad ;\quad \Im \rho_k =\CO(e^{-2(S_0-\eta)/h}),
\ee
uniformly.
In particular, taking $\mu$ and $\widetilde \mu$ as in Theorem \ref{Th5}, the result easily follows. Moreover, since the previous discussion can be applied to any $\mu'\in [e^{-\eta /h}, h^\delta]$,  application of Theorem \ref{Th4} tells us that the resonances of $P$ in $[\lambda_0, \lambda_0+ Ch]-i[0, \frac12 h^{2n+\max(\frac{n}2,1) +1+3\delta}]$ satisfy to the same estimates (\ref{estresF}).
 
  \appendix
%%%%%%%%%%%%%%%%%%%%%%%%%%%  APPENDIX %%%%%%%%%%%%%%%%%%%%%%%%%%%
\section{Appendix}\label{sec-app}
\subsection{Proof of Lemma 5.1}\label{app1}
We denote by $\chi_0$ a  real smooth function on $\R$ verifying,
\begin{itemize}
\item $\chi_0 (s) =0$ for $s\leq 0$;
\item $\chi_0(s) = 1$ for $s\geq \ln 2$;
\item $\chi_0$ is non decreasing.
\end{itemize}
Then, for $r\geq 0$, we set,
$$
G(r) := \chi_0(r-R_0)(1-\chi_0(r-\ln\lambda ))e^r + 2\lambda \chi_0(r-\ln\lambda ),
$$
and,
$$
g(r):= \int_0^r G(s)ds.
$$
In particular, $g$ verifies Condition (i) of Lemma \ref{flambda}, and we have,
\begin{itemize}
\item $G(r) =\chi_0(r-R_0)e^r$ for $r\in [R_0, \ln\lambda ]$;
\item $G(r) =(1-\chi_0(r-\ln\lambda ))e^r + 2\lambda \chi_0(r-\ln\lambda )$ for $r\in [\ln\lambda ,\ln 2\lambda] $;
\item $G(r) =2\lambda $ for $r\in [\ln 2\lambda, +\infty)$.
\end{itemize}
Thus, $g'=G\leq 2\lambda$  and $g''(r)=G'(r)\geq 0$ on $\R_+$ (this is immediate on $[R_0, \ln\lambda ]\cup [\ln 2\lambda, +\infty)$, while, on $[\ln\lambda ,\ln 2\lambda]$, we compute,
$G'(r) = (1-\chi_0(r-\ln\lambda ))e^r + \chi_0'(r-\ln\lambda )(2\lambda -e^r)\geq 0$). 

Therefore, $g$ is convex on $\R_+$, so that Condition (iii) of Lemma \ref{flambda} is verified by $g$, too, while Condition (v) is obvious. 

As for condition (iv), we observe,
\begin{itemize}
\item On $[0, R_0+\ln 2]$, one has, $g'+|g''| = \CO(1)$;
\item On $[R_0+\ln 2, \ln\lambda]$, one has,
$g(r)\geq \int_{R_0+\ln 2}^r e^sds = e^r -2e^{R_0}$, while $g'(r)=g''(r) =e^r\leq g(r)+2e^{R_0}$;
\item On $[\ln\lambda , +\infty)$, one has, $g(r)\geq g(\ln\lambda)=\lambda$, and thus  $g'+|g''| = \CO(g)$.
\end{itemize}
So,  $g$ verifies Conditions (ii)-(v) of Lemma \ref{flambda}. 

For $r\in [\ln2\lambda, +\infty)$, we have,
\be
\label{gr+}
g(r) = g(\ln 2\lambda) + 2\lambda (r -\ln2\lambda) = 2\lambda r - \alpha_\lambda,
\ee
where $\alpha_\lambda:=  2\lambda\ln 2 \lambda -g(\ln 2\lambda)$,
and, since,
\begin{eqnarray*}
 && g(\ln 2\lambda)\leq \int_0^{\ln\lambda} e^r dr + \int_{\ln\lambda}^{\ln 2\lambda} 2\lambda dr = (1+2\ln 2)\lambda;\\
 && g(\ln 2\lambda)\geq \int_{R_0+\ln 2}^{\ln 2\lambda}e^rdr \geq 2\lambda -2e^{R_0}.
\end{eqnarray*}
we see that,
$$
2\lambda\ln 2 \lambda -(1+2\ln 2)\lambda\leq\alpha_\lambda \leq 2\lambda\ln 2 \lambda - 2 \lambda +2e^{R_0}.
$$
Therefore, for $\lambda$ large enough, the unique point $r_\lambda$, solution of $g(r_\lambda) = \lambda r_\lambda$, is given by,
\be
\label{rlambda}
r_\lambda =\frac{\alpha_\lambda}{\lambda} \in [2\ln \lambda -1, 2\ln  \lambda -  2+2\ln 2+2\lambda^{-1}e^{R_0}]\subset [2\ln\lambda -1, 2\ln\lambda -\varepsilon_0],
\ee
where $\varepsilon_0 := 1-\ln2 >0$.

Now, we fix some real-valued function $\varphi_0\in C^\infty (\R)$, such that,
\begin{itemize}
\item $\varphi_0 (s) = 2s$ for $s\leq -\varepsilon_0$;
\item $\varphi_0 (s) = s$ for $s\geq \varepsilon_0$;
\item $1\leq \varphi_0'\leq 2$ everywhere.
\end{itemize}
Then, using (\ref{gr+})-(\ref{rlambda}), we see  that the function $f_\lambda$ defined by,
\begin{itemize}
\item $f_\lambda (r):= g(r)$ for $r\in [0, \ln 2\lambda]$;
\item $f_\lambda (r):= \lambda \varphi_0 (r-r_\lambda) + \alpha_\lambda$ for $r\geq \ln2\lambda$,
\end{itemize}
is smooth on $\R_+$, and verifies all the conditions required in Lemma \ref{flambda}. \hfill$\square$

\subsection{The distorted Laplacian}
\begin{lem}\sl
\label{lemA1}
If $\theta>0$ is small enough, the function $\Phi_{\theta}$ defined in (\ref{defdist}) verifies,
$$
\Im [({}^td\Phi_{\theta}(x))^{-1}\xi]^2 \leq -\theta a (|x|)|\xi|^2.
$$
for all $(x,\xi)\in\R^{2n}$.
\end{lem}
\begin{proof} Let $F:= {}^tdA= dA =(F_{i,j})_{1\leq i,j\leq n}$. We compute,
$$
F_{i,j}(x)=a(x)\delta_{i,j} +a' (|x|)\frac{x_ix_j}{|x|},
$$
that is, denoting by $\pi_x :=|x|^{-2} x\cdot{}^tx$ the orthogonal projection onto $\R x$, and recalling the notation $b(r)= r a (r)$,
$$
F (x)=a(|x|)I + a'(|x|)|x|\pi_x = b' (|x|)\pi_x + a(|x|)(I-\pi_x).
$$
In particular, using Lemma \ref{flambda}, we obtain,
$$
0\leq a(|x|) \leq F (x) \leq 2,
$$
in the sense of self-adjoint matrices.
On the other hand, we have,
$$
({}^td\Phi_{\theta}(x))^2 =(I+i \theta F(x))^2 = S_\theta+i T_\theta,
$$
with $S_\theta=I -\theta^2F(x)^2$ and $T_\theta=2 \theta F(x)$. Hence,  $T_\theta\geq 0$, and, since $S_\theta$, $T_\theta$ and $F$ commute, an easy computation gives,
\begin{eqnarray*}
\Im [({}^td\Phi_{\theta}(x))^{-1}\xi]^2
&=& -T_\theta(S_\theta^2+T_\theta^2)^{-1}\xi\cdot\xi
= - 2\theta F (1+\theta^2F^2)^{-2}\xi\cdot\xi .
\end{eqnarray*}
As a consequence, for $\theta$ small enough, we find,
$$
\Im [({}^td\Phi_{\theta}(x))^{-1}\xi]^2\leq -\theta F(x)\xi\cdot\xi\leq -\theta a(|x|)|\xi|^2.
$$
\end{proof}

%%%%%%%%%%%%%%%%%%  BIBLIOGRAPHY  %%%%%%%%%%%%%%%%%%%%%%%%%%%%%%%%

%\bibliographystyle{amsplain}

%\bibliography{mrs}

\begin{thebibliography}{10}

\bibitem{AgCo}
J.~Aguilar and Jean-Michel Combes, \emph{A class of analytic perturbations for
  one-body Schr\"odinger Hamiltonians}, Comm. Math. Phys. 22 (1971), 269--279.

\bibitem{BaCo}
Erik Balslev and Jean-Michel Combes, \emph{Spectral properties of many-body
  Schr\"odinger operators with dilation analytic interactions}, Comm. Math.
  Phys. 22 (1971), 280--294.

\bibitem{CaMaRa}
Claudy Cancelier, Andr{\'e} Martinez, and Thierry Ramond, \emph{Quantum
  resonances without analyticity}, Asymptotic Analysis. \textbf{44} (2005),
  no.~1-2, 47--74.

\bibitem{Cy}
Hans~L. Cycon, \emph{Resonances defined by modified dilations}, Helv. Phys.
  Acta 58, 969--981 (1985).

\bibitem{DiSj}
Mouez Dimassi and Johannes Sj{\"o}strand, \emph{Spectral asymptotics in the
  semi-classical limit}, London Mathematical Society Lecture Note Series, vol.
  268, Cambridge University Press, Cambridge, 1999.

\bibitem{FuLaMa}
Setsuro Fujii{\'e}, Amina Lahmar-Benbernou, and Andr{\'e} Martinez, \emph{Width
  of shape resonances for non globally analytic potentials}, in preparation.

\bibitem{ga}
George Gamow, \emph{{Zur Quantentheorie der Atomzertr\"ummerung.}}, Z. f.
  Physik \textbf{52} (1928), 510--515.

\bibitem{GeMa}
Christian G{\'e}rard and Andr{\'e} Martinez, \emph{Principe d'absorption limite
  pour des op\'erateurs de {S}chr\"odinger \`a longue port\'ee}, C. R. Acad.
  Sci. Paris S\'er. I Math. \textbf{306} (1988), no.~3, 121--123.

\bibitem{GeSi}
Christian G\'erard and Israel~M. Sigal, \emph{Space-time picture of semiclassical
  resonances}, Comm. Math. Phys. 145 (1992), 281--328.

\bibitem{HeSj}
Bernard Helffer and Johannes Sj{\"o}strand, \emph{R{\'e}sonances en limite
  semi-classique}, M{\'e}m. Soc. Math. France (N.S.) \textbf{24-25} (1986),
  iv+228.

\bibitem{Hu}
Walter Hunziker, \emph{Distortion analyticity and molecular resonance curves},
  Ann. Inst. H. Poincar\'e Phys. Th\'eor. \textbf{45} (1986), no.~4, 339--358.

\bibitem{JeNe}
Arne Jensen and Gheorghe Nenciu, \emph{The {F}ermi golden rule and its form at
  thresholds in odd dimensions}, Comm. Math. Phys. \textbf{261} (2006), no.~3,
  693--727. 

\bibitem{Ma2}
Andr{\'e} Martinez, \emph{Resonance free domains for non globally analytic
  potentials}, Ann. Henri Poincar\'e 3 (2002), no. 4, 739--756 ; Erratum 8
  (2007), 1425-1431.

\bibitem{Ma1}
\bysame, \emph{An introduction to semiclassical and microlocal analysis},
  Universitext, Springer-Verlag, New York, 2002.

\bibitem{MeSj}
Anders Melin and Johannes Sj{\"o}strand, \emph{Fourier integral operators with
  complex-valued phase functions}, Fourier integral operators and partial
  differential equations (Colloq. Internat., Univ. Nice, Nice, 1974), Springer,
  Berlin, 1975, pp.~120--223. Lecture Notes in Math., Vol. 459.

\bibitem{na1}
Shu Nakamura, \emph{Distorsion analyticity for two-body Schr\"odinger
  operators}, Ann. Inst. H. Poincar\'e, Phys. Th\'eor. 53, 149-157 (1990).

\bibitem{na2}
\bysame, \emph{Shape resonances for distortion analytic Schr\"odinger
  operators}, CPDE 14(10), 1989, 1385--1419.

\bibitem{Or}
Andreas Orth, \emph{Quantum mechanical resonance and limiting absorption: the
  many body problem}, Comm. Math. Phys. \textbf{126} (1990), no.~3, 559--573.

\bibitem{Sig}
Israel~M. Sigal, \emph{Complex transformation method and resonances in one-body
  quantum systems}, Ann. Inst. H. Poincar\'e, Phys. Th\'eor. 41, 103--114 (1984);
  Addendum 41, 333 (1984).

\bibitem{Sim}
Barry Simon, \emph{The definition of molecular resonance curves by the method
  of exterior complex scaling}, Phys. Letts. 71A, 211--214 (1979).

\bibitem{SoWe}
Avy Soffer and Michael~I. Weinstein, \emph{{Time dependent resonance theory.}},
  Geom. Funct. Anal. \textbf{8} (1998), no.~6, 1086--1128.

\bibitem{TaZw}
Siu-Hung Tang and Maciej Zworski, \emph{From quasimodes to resonances}, Math.
  Res. Lett. \textbf{5} (1998), no.~3, 261--272.

\end{thebibliography}

%Copie de  MRS12.bbl insereree le 9 mai a 8h40%%%%%%%%%%%%%

\providecommand{\bysame}{\leavevmode\hbox to3em{\hrulefill}\thinspace}
\providecommand{\MR}{\relax\ifhmode\unskip\space\fi MR }
\providecommand{\MRhref}[2]{%
  \href{http://www.ams.org/mathscinet-getitem?mr=#1}{#2}
}
\providecommand{\href}[2]{#2}

\end{document}